\baselineskip=15pt plus 2pt
\magnification=1100
\font\medtenrm=cmr10 scaled \magstep2

\def\sqr#1#2{{\vcenter{\vbox{\hrule height.#2pt\hbox{\vrule width.#2pt
height#1pt\kern#1pt \vrule width.#2pt}\hrule height.#2pt}}}}
\def\square{\mathchoice\sqr64\sqr64\sqr{2.1}3\sqr{1.5}3}

\centerline {\medtenrm On linear systems and $\tau$ functions 
associated with Lam\'e's equation}\par
\centerline {\medtenrm and Painlev\'e's equation VI}\par
\vskip.05in
\centerline {\medtenrm Gordon Blower}\par
\vskip.05in
\centerline {\sl Department of Mathematics and Statistics,}\par
\centerline {\sl Lancaster
University, Lancaster LA1 4YF, England UK}\par
\vskip.05in

\centerline {17th November 2009}\par
\vskip.1in
\hrule
\vskip.1in
\noindent {ABSTRACT} Painlev\'e's transcendental differential
equation $P_{VI}$ may be expressed as the consistency condition
for a pair of linear differential equations with $2\times 2$ matrix
coefficients with rational entries. By a construction due to Tracy
and Widom, this linear system is associated
with certain kernels which give trace class operators on Hilbert space.
 This
paper expresses such operators in terms of Hankel operators
$\Gamma_\phi$ of linear 
systems which are realised in terms of the Laurent coefficients of
the solutions of the differential equations. Let $P_{(t, \infty
)}:L^2(0, \infty )\rightarrow L^2(t, \infty )$ be the orthogonal
projection. For such, the Fredholm determinant $\tau (t)=\det
(I-P_{(t, \infty )}\Gamma_\phi
)$ defines the $\tau$
function, which is here expressed in terms of the solution of a matrix Gelfand Levitan
equation. For suitable values of the parameters, solutions of the
hypergeometric equation give a linear system with similar
properties. For meromorphic transfer functions $\hat \phi$ that have
poles on an arithmetic progression, the corresponding Hankel
operator has a simple form with respect to an exponential basis in
$L^2(0, \infty )$; so $\det (I-\Gamma_\phi P_{(t, \infty )})$ can be expressed as a
series of finite determinants. This applies to elliptic functions of
the second kind, such as satisfy Lam\'e's equation with $\ell =1$.\par
\par
\vskip.05in
\noindent MSC2000 classification 47B35\par

\vskip.1in
\hrule
\vskip.05in
\noindent {\bf 1. Introduction}\par
\vskip.05in

\indent Tracy and Widom [29] observed that many important kernels in
random matrix theory arise from solutions of linear differential equations
with rational coefficients. In particular, the classical systems of
orthogonal polynomials can be expressed in such terms. In this paper, we extend the scope of
their investigation by analysing kernels associated with Lam\'e's
equation and Painlev\'e's equation VI. As these differential equations
have solutions which may be expressed in terms of elliptic functions,
we begin by reviewing and extending the definitions from [29].\par
\indent Let $P(x,y)$ be an irreducible complex
polynomial, and $n$ the degree of $P(x,y)$ as a polynomial in $y$. 
Then we introduce the curve
${\cal E}=\{ (\lambda , \mu )\in {\bf C}:P(\lambda , \mu )=0\}$, and
observe that ${\cal E}\cup\{ (\infty ,\infty )\}$ gives a compact
Riemann surface which is the $n$-sheeted branched cover of
Riemann's sphere ${\bf
P}^1$. Let
${\bf K}$ be splitting field of $P(x,y)$ over ${\bf C}(x)$, so we
can regard ${\bf K}$ as the space of functions of rational character 
on ${\cal
E}.$ Let $g$ be the genus of ${\cal E}$, and introduce the Jacobi variety 
${\bf J}$ of ${\cal E}$,
which is the quotient of ${\bf C}^g$ by some lattice ${\bf L}$ in
${\bf C}^g$.\par
\vskip.05in
\noindent {\bf Definition.} By a Tracy--Widom system [29] we mean a differential equation
$${{d}\over{dx}}\left[ \matrix{ f\cr g\cr}\right] =\left[ \matrix
{\alpha &\beta \cr -\gamma & -\alpha\cr}\right] \left[ \matrix{
f\cr g\cr}\right]\eqno(1.1)$$
\noindent where $\alpha, \beta, \gamma$ belong to ${\bf K}$ or more
generally are locally rational
functions on ${\bf J}$. Then for solutions with
 $f,g\in L^\infty ((0, \infty ); {\bf R})$, we introduce an integrable operator on $L^2(0, \infty )$
by the kernel
$$K(x,y)={{f(x)g(y)-f(y)g(x)}\over {x-y}}\qquad (x\neq y;
x,y\in {\bf R})\eqno(1.2)$$
\indent The kernel $K$ compresses to give an integral operator $K_S$ on $L^2(S; dx)$
for any subinterval $S$ of $(0, \infty )$ and it is important to identify
those $K_S$ such that $K_S$ is of trace class and $0\leq K_S\leq I$. In
such cases, the Fredholm determinant $\det (I+\lambda K_S)$ is defined
and $K_S$ is associated with a determinantal random point field on $S$.
In particular, $\det (I-K_{(t, \infty )})$ gives the probability that
there are no random points on $(0, \infty ).$ \par
\vskip.05in
\noindent {\bf Definition} {\sl ($\tau$-function).} Suppose that
$K:L^2(0,\infty )\rightarrow L^2(0, \infty )$ is a self-adjoint operator
such that $K\leq I$, $K$ is trace class and $I-K$ is invertible.
For a measurable subset $S$ of $(0,\infty )$, let $P_S:L^2(0, \infty
)\rightarrow L^2(S)$ be the orthogonal projection given by $f\mapsto f{\bf I}_S$, where ${\bf I}_S$
is the indicator function of $S$. Then the $\tau$ function is
$$\tau (t)=\det (I-KP_{[t, \infty )})\qquad (t>0).\eqno(1.3)$$

\vskip.05in
\indent The purpose of this paper is to take 
kernels that are given by certain Tracy--Widom systems, and show how
to express the corresponding $\tau$ in terms of the solution of a
Gelfand--Levitan integral equation. Our technique involves linear systems, and extends ideas
developed in [6], and leads to a solution of the integral equation
in terms of the linear system.\par
\indent Let $H$ be a complex separable Hilbert spaces, known as the
state space, and
let $(e^{-tA})_{t>0}$ a bounded $C_0$-semigroup of linear operators on $H$;
so that ${\cal D}(A)$ is a dense linear subspace of $H$, and $\Vert
e^{-tA}\Vert \leq M$ for all $t>0$ and some $M<\infty$. Then
let $B:{\bf C}\rightarrow {\cal D}(A)$ and $C:{\cal D}(A)\rightarrow
{\bf C}$ be bounded linear operators, and introduce the linear system
$$\eqalignno{ {{dX}\over {dx}}&=-AX+BU\qquad (X(0)=0),\cr
                            Y&=CX&(1.4)\cr}$$

\noindent known as $(-A,B,C)$. Under further conditions to be discussed below, the integral
$$R_x=\int_x^\infty e^{-tA}BCe^{-tA}\,dt\eqno(1.5)$$
\noindent converges and defines a trace class operator on $H$. The
notation suggests that $R_x$ is a resolvent operator.\par

\vskip.05in
\noindent {\bf Definition} {\sl (Hankel operator).} For a linear system as above, we introduce the {\sl symbol} 
$\phi (x)=Ce^{-xA}B$, which gives a
bounded function $\phi :(0, \infty )\rightarrow {\bf C}$; this term should not be
confused with the different usage in [26, p 6]. Generally, for
$E$ a separable complex Hilbert space and $\phi\in L^{2}((0, \infty
);E)$, let $\Gamma_\phi$ be the
Hankel operator 
$$\Gamma_\phi h(x)=\int_0^\infty \phi (x+y)h(y)\, dy\eqno(1.6)$$
\noindent defined on a suitable domain in $L^2(0, \infty )$ into
$L^2((0, \infty );E)$. \par
\indent By forming orthogonal sums of the state space and block
operators, we can form sums of symbol functions. Likewise, by forming
tensor products of state spaces and operators, we can from products
of symbol function. Using these two basic constructions, we can form
some apparently complicated symbol functions, starting from the basic
multiplication operator $A: f(t)\mapsto tf(t)$ in $L^2(0, \infty )$.
Thus we extend the method of section 2 to a more intricate problem.\par
\par
\indent In section 3, we consider operators related to the 
solution of
Painlev\'e's transcendental equation VI
$$\eqalignno{{{d^2y}\over{dt^2}}&={{1}\over{2}}\Bigl(
{{1}\over{y}}+{{1}\over{y-1}}+{{1}\over{y-t}}\Bigr)\Bigl({{dy}\over{dt}}
\Bigr)^2-\Bigl( {{1}\over{t}}+{{1}\over{t-1}}+{{1}\over{y-t}}\Bigr)
{{dy}\over{dt}}\cr
&\quad + {{y(y-1)(y-t)}\over{t^2(t-1)^2}}\Bigl( \alpha +{{\beta
t}\over{y^2}}+{{\gamma (t-1)}\over{(y-1)^2}}+{{\delta
t(t-1)}\over{(y-t)^2}}\Bigr).&(1.7)\cr}$$
\noindent with constants
$$\alpha ={{1}\over{2}}(\theta_\infty -1)^2,\quad \beta
=-{{1}\over{2}}\theta_0^2,\quad \gamma ={{1}\over{2}}\theta_1^2,\qquad
\delta ={{1}\over{2}}(1-\theta_t^2)\eqno(1.8)$$
\noindent and 
$$\theta_\infty
=-2(z_0+z_1+z_t)-(\theta_0+\theta_1+\theta_t).\eqno(1.9)$$

\noindent Jimbo, Miwa and Ueno [15, 16] showed that 
the nonlinear differential equation $P_{VI}$ is the compatibility
condition for the pair of linear differential equations
$$\eqalignno{{{d\Phi}\over{d\lambda}}&=\Bigl({{W_0}\over{\lambda}}+
{{W_1}
\over{\lambda -1}}+{{W_t}\over{\lambda -t}}\Bigr)\Phi&(1.10)\cr
{{d\Phi}\over{dt}}&={{-W_t}\over{\lambda -t}}\Phi&(1.11)\cr}$$
\noindent on the punctured Riemann sphere with $2\times 2$ complex matrices 
$W_0, W_1, W_t$ depending upon $t$; see (3.8) for the entries. Using the Laurent series 
of $\Phi (\lambda )$, we introduce a linear system $(-A,B,C)$ that realises
$\Phi$ and deduce information about the Hankel operator
$\Gamma_\Phi$. In previous papers [5,6], we have considered kernels that
factorize as $K=\Gamma_\phi^\dagger\Gamma_\phi$ where $\Gamma_\phi$ is
Hilbert-Schmidt, so that $K\geq 0$ and $K$ is trace class. In the
context of $P_{VI}$, we show that the prescription (1.2) gives a kernel
$K$ that admits a factorization
$K=\Gamma_\phi^\dagger \sigma\Gamma_\phi$, where $\sigma$ is a constant
signature matrix. In section 5 we introduce a suitable
$\tau$ function and express this in terms of the solution of an integral
equation of Gelfand--Levitan type, which we can solve in terms of
the linear system. A similar approach works for suitable solutions of 
Gauss's hypergeometric equation with a restricted choice of
parameters, as we show in section 5.\par

\vskip.05in
\noindent {\bf Definition} {\sl (Transfer function).} Given a 
Hilbert space $E$, for $\phi\in L^2((0, \infty ); 
dt; E)$ let 
$$\hat \phi (s)=\int_0^\infty e^{-st}\phi (t)\, dt\eqno(1.12)$$
\noindent be the transfer function of $\phi$, otherwise known as
the Laplace transform, which gives an analytic
function from $\{ s:\Re  s>0\}$ into $E$.\par
\indent We assume that $\hat\phi$ is meromorphic, and that, by virtue 
of the Mittag-Leffler theorem, one can express $\phi$
as a series
$$\phi (x)=\sum_{j=1}^\infty \xi_je^{-\lambda_jx}\eqno(1.13)$$
\noindent in which we shall always assume that $\Re \lambda_j>0$ and that the
$e^{-\lambda_jx}$ are linearly independent in $L^2(0, \infty )$. We
wish to express various $\tau$ functions in terms of the determinants 
$$D_{S\times T}=\det \Bigl[ {{1}\over{\lambda_j+\bar\lambda_k}}
\Bigr]_{(j,k)\in S\times T}\eqno(1.14)$$
\noindent where $S$ and $T$ are finite subsets of ${\bf N}$ of
equal cardinality. In sections 6, we consider Hankel operators
with symbols as in (1.13), and establish basic results about the 
expansions of $\det (I-\Gamma_\phi )$ in terms of the bases. In
particular, if $(\lambda_j)_{j=1}^\infty$ forms an arithmetic progression in the
plane, then 
$\hat\phi (s)=\sum_{j=1}^\infty {{\xi_j}\over{s+\lambda_j}}$ gives a cardinal series.\par 
\indent In section 7, we consider the Bessel kernel, which arises in 
random matrix theory as the hard edge of the eigenvalue distribution
from the Jacobi ensemble [28]. Let $J_\nu$ be
Bessel's function of the first kind of order $\nu$, and let
$u(x)=\sqrt{x}J_\nu (2\sqrt{x})$, which satisfies
$${{d^2u}\over{dx^2}}+\Bigl(
{{1}\over{x}}+{{1-\nu^2}\over{4x^2}}\Bigr)u(x)=0.\eqno(1.15)$$
\noindent We introduce $\phi (x)=u(e^{-x})$, and the Hankel operator
$\Gamma_\phi$ with symbol $\phi$. The transfer function $\hat \phi$
is meromorphic with poles on an arithmetic progression on the
positive real axis, so we are able to obtain a simple expansion for 
$\tau (t)=\det (I-\Gamma_\phi^2P_{[t, \infty )}),$ and identify the
determinants $D_{N\times N}$ with combinatorial objects.\par
\indent In
section 8 we consider solutions of Lam\'e's equation 
$$\Bigl( -{{d^2}\over{dz^2}}+\ell (\ell
+1)k^2{\hbox{sn}}(z\mid k)^2\Bigr) \Phi (z)=\lambda \Phi
(z)\eqno(1.16)$$ 
\noindent which we express as a differential equation on 
the elliptic curve
$Z^2=4(X-e_1)(X-e_2)(X-e_3)$. The solution gives rise to an elliptic
function $\phi$ such that $\hat \phi$ has poles on a bilateral
arithmetic procession parallel to the imaginary axis in ${\bf C}$. Hence we can prove
results concerning the Fredholm determinant of $\Gamma_\phi$.\par

\vskip.05in
\vskip.05in

\noindent {\bf 2. The $\tau$ function associated with a linear
system}\par
\vskip.05in
\noindent In this section we introduce the basic example of the linear
system which we will use in sections 3 and 5 to realise solutions of some
differential equations. In [30], Tracy and Widom consider physical
applications of the kernels $R_x$ that we introduce her.\par
\vskip.05in

\noindent {\bf Definition} {\sl (Integrable operators)}. Let $f_1,
\dots , f_N, g_1, \dots ,g_N\in L^\infty (0, \infty )$ satisfy 
$$\sum_{j=1}^N f_j(x)g_j(x)=0\qquad (x>0).$$
\noindent Then the integral operator $K$ on $L^2(0, \infty )$ that has 
kernel 
$$K\leftrightarrow {{\sum_{j=1}^Nf_j(x)g_j(y)}\over{x-y}}\eqno(2.1)$$
\noindent is said to be an integrable operator; see [9]. One can show that $K$
is bounded on $L^2(0, \infty )$.\par
\indent Let ${\cal D}(A)=\{ f\in L^2(0, \infty ): tf(t)\in L^2(0, \infty )\}$
and for $b,c\in {\cal D}(A)$ introduce the operators:
$$\matrix{ A:&{\cal D}\subset L^2(0, \infty )\rightarrow 
L^2(0, \infty ):&f(x)\mapsto xf(x)\cr
B:&{\bf C}\rightarrow {\cal D}(A):& \alpha\mapsto b\alpha;\cr
  C:&{\cal D}(A)\rightarrow {\bf C}: 
&f\mapsto \int_0^\infty f(s)c(s)\, ds\cr
\Theta_x:&L^2(0, \infty )\rightarrow L^2(0,\infty ) :& \Theta_xf(t)=e^{-xt}\bar
c(t)\hat f(t)\cr
\Xi_x:&L^2(0, \infty )\rightarrow L^2(0, \infty ): & \Xi_x
f(t)=e^{-xt}b(t)\hat f(s)\cr}\eqno(2.2)$$ 
\noindent Then we introduce $\phi (s)=Ce^{-sA}B$ and
$\phi_{(x)}(s)=\phi (s+2x)$, and the Hankel integral operator
$\Gamma_{\phi_{(x)}} $ with kernel $\phi (s+t+2x).$ Then we introduce 
$R_x=\int_x^\infty e^{-tA}BCe^{-tA}\, dt$ which has kernel
$$ R_x\leftrightarrow {{b(t)c(s)e^{-x(s+t)}}\over{s+t}}\qquad
(s,t>0).\eqno(2.3)$$
\vskip.05in
\noindent {\bf Proposition 2.1.} {\sl Suppose that $c(t)/\sqrt {t}$ and
$b(t)/\sqrt{t}$ belong to $L^2(0, \infty )$, and that $c$ and $b$
belong to $L^\infty (0, \infty )$.\par
\indent (i) Then
$\Gamma_{\phi_{(x)}}$ and $R_x$ are trace class operators for
all $x\geq 0$.\par
\indent (ii) Suppose further that $I+\lambda R_x$ is
invertible for some $\lambda \in {\bf C}$. Then 
$$T_\lambda (x,y)=
-\lambda Ce^{-xA}(I+\lambda R_x)^{-1}e^{-yA}B\eqno(2.4)$$
\noindent gives the solution to the equation
$$\lambda \phi (x+y)+T_\lambda (x,y)+\lambda \int_x^\infty T_\lambda
(x,z)\phi (z+y)\, dz=0\qquad (0<x<y)\eqno(2.5)$$
\noindent and}
$$T_\lambda (x,x)={{d}\over{dx}}\log\det (I+\lambda
\Gamma_{\phi_{(x)}}).\eqno(2.6)$$
\indent {\sl (iii) The operator $R_x^2$ is an integrable operator with kernel
$$R_x^2\leftrightarrow e^{-xu}b(u)
{{f_x(u)-f_x(t)}\over
{t-u}}c(t)e^{-xt}\eqno(2.7)$$
\noindent where} 
$$f_x(u)=\int_0^\infty {{b(t)c(t)e^{-tx}}\over {u+t}}\, dt.\eqno(2.8)$$
\indent {\sl (iv) If $I+\lambda R_x$ and $I-\lambda R_x$ are invertible, 
then there exists an integrable operator $L_x(\lambda )$ such that}
$$I+L_x(\lambda )=(I-\lambda^2R_x^2)^{-1}.\eqno(2.9)$$
\vskip.05in
\noindent {\bf Proof.} (i) One checks that $\Theta_x$ has kernel
$e^{-st}e^{-xt}\bar c(t)$ and that $\Xi_x$ has kernel
$e^{-st-xs}b(s)$; hence $\Theta^\dagger_x$ and $\Xi_x$ are
Hilbert--Schmidt operators. One verifies that their products are 
$R_x=\Xi_x\Theta_x^\dagger$ and
$\Gamma_{\phi_{x}} =\Theta^\dagger\Xi_x$, and hence $R_x$ and $\Gamma_x$
are trace class.\par
\indent (ii) Using (i), we can check that  $\det (I+\lambda R_x)=\det
(I+\lambda \Gamma_{\phi_{(x)}})$. Then one verifies the remainder by
using Lemma 5.1(iii) of [6].\par 
\indent (iii) This result is essentially contained in lemma 2.18
of [9], but we give a proof for completeness. The kernel of $R_x^2$ is
$$b(s)e^{-sx}c(u)e^{-ux}\int_0^\infty {{b(t) c(t)
e^{-2tx}}\over {(s+t)(u+t)}} dt\qquad (u,s>0),\eqno(2.10)$$
\noindent and one can decompose this expression by using partial
fractions. By the Cauchy--Schwarz inequality, $\vert f_x(u)\vert^2 \leq
\int_0^\infty t^{-1}b(t)^2dt\int_0^\infty t^{-1}c(t)^2dt$, so $f_x$ is
bounded.\par 
\indent (iv) Furthermore, $(I-\lambda R_x)^{-1}(I+\lambda
R_x)^{-1}$ is a bounded linear operator; so by Lemma 2.8 of [9], there exists an integrable
operator $L_x$ such that $(I+L_x(\lambda ))(I-\lambda^2R_x^2)=I.$\par
\rightline{$\square$}\par
\noindent Given an integrable operator $K$ on $L^2(a,b)$ such that
$I-K$ is invertible, the authors of [9] show how to express
$(I-K)^{-1}$ as the solution of a Riemann--Hilbert problem on a bounded interval


\vskip.05in
\noindent {\bf 3. A linear system associated with Painlev\'e's 
equation VI}\par
\vskip.05in
\indent The Painlev\'e equation $P_{VI}$ is associated with the system
$$\eqalignno{{{d\Phi}\over{d\lambda}}&=\Bigl({{W_0}\over{\lambda}}+
{{W_1}
\over{\lambda -1}}+{{W_t}\over{\lambda -t}}\Bigr)\Phi &(3.1)\cr
{{d\Phi}\over{dt}}&={{-W_t}\over{\lambda -t}}\Phi&(3.2)\cr}$$
\noindent where the fixed singular points are $\{0 , 1,\infty \}$ and 
$$W_\nu =W_\nu (t)=\left[\matrix{ z_\nu+\theta_\nu /2 &
- u_\nu z_\nu\cr u_\nu^{-1}(z_\nu
+\theta_\nu)& -z_\nu-\theta_\nu /2\cr}\right]\qquad (\nu
=0,1,t)\eqno(3.3)$$
\noindent with parameters $\theta_\nu$ and $z_\nu$ satisfying various
conditions specified in [16]. The
consistency condition for the system (3.1) and (3.2) reduces to the identity 
$${{1}\over{\lambda}}{{\partial W_0}\over{\partial
t}}+{{1}\over{(\lambda -1)}}{{\partial W_1}\over{\partial
t}}+{{1}\over{(\lambda -t)}}{{\partial W_t}\over{\partial
t}}={{[W_0,W_t]}\over{\lambda
(\lambda -t)}}+{{[W_1,W_t]}\over{(\lambda -1)(\lambda -t)}},\eqno(3.4)$$
\noindent which leads, after a lengthy computation given in
Appendix C of [16], to the equation $P_{VI}$.\par 
\indent Jimbo {\sl et al} [15, 16, 17] introduced pairs of differential 
equations (3.1) and (3.2) such that (3.5) reduces to one of the Painlev\'e equations. In
the present context (3.1) are known as the deformation equations
and (3.4) is associated with the names of Schlesinger and
Garnier [10]. Note that ${\hbox {trace}}\, W=0 $ if and only if $JW$ is
symmetric; also $W$ is nilpotent if and only if $JW$ is
symmetric and $\det (JW)=0$.\par

\indent First we introduce a linear system for the differential
equation (3.4); later we introduce a linear system that realises the
kernel most naturally associated with $P_{VI}$. For notational
simplicity, we often suppress the dependence of operators upon $t$. The following result 
is a consequence of results of Turrittin [31, 27],
who clarified certain facts about the Birkhoff canonical form for
matrices.\par
\vskip.05in
\noindent {\bf Lemma 3.1.} {\sl Let $W_\infty =-(W_0+W_1+W_t)$
and suppose that the eigenvalues of
$W_\infty$ are $\mp \theta_\infty/2$ where $\pm \theta_\infty$ is 
not a positive integer, and let $\Phi_0$ be a constant $2\times 1$
vector. Then there exist 
$2\times 2$ complex matrices $C_j$ for $j=1, 2,\dots $, depending upon $t$, such that
$$\Phi (x)=\Bigl( I+\sum_{j=1}^\infty {{C_j}\over{x^j}}\Bigr)
x^{-W_\infty}\Phi_0\qquad (\vert x\vert >t )\eqno(3.5)$$
\noindent satisfies the differential equation (3.1).}\par
\vskip.05in
\noindent {\bf Proof.} We can define 
$x^{-W_\infty}=\exp (-W_\infty\log x)$ as a convergent power
series. By considering terms in 
the convergent Laurent series,
one requires to show that there exist coefficients $C_0=I$ and $C_j$ that satisfy the recurrence relation
$$\eqalignno{C_n(-W_\infty-nI)&=-W_\infty C_n+W_1(C_0+\dots
+C_{n-1})\cr
&\quad +tW_t(t^{n-1}C_0+t^{n-2}C_1+\dots +C_{n-1}),&(3.6)}$$
\noindent where $W_\infty +nI$ and $W_\infty$ have no common
eigenvalues. Sylvester showed that, given square matrices $V, W$ and
$Z$ such that $V$ and $W$
have no eigenvalues in common, the matrix equation $CV-WC=Z$ has a
unique solution $C$; see [31, Lemma 1]. Hence unique $C_n$ exist,
and one shows by induction
that $\Vert C_n\Vert$ is at most of geometric growth in $n$. In
particular, if $\Vert W_\infty \Vert <1$, then the solution of
$W_\infty C_n-C_n(W_\infty +nI)=D_n$ is
$$C_n=-\int_0^\infty e^{sW_\infty }D_ne^{-s(W_\infty +nI)}\,
ds.\eqno(3.7)$$
\rightline{$\square$}\par
\par 
\vskip.05in
\indent We have proved that (3.1) has a solution in a neighbourhood
of infinity, and one can show that it extends to an analytic
solution on the universal cover of the punctured Riemann sphere
${\bf P}^1\setminus \{ 0,1,t, \infty \}$. (Jimbo, Miwa and Ueno [15] have shown that any $C^2$ solution of
the pair (3.1) and (3.2) on ${\bf R}$ extends to a meromorphic solution on ${\bf
C}$; see [10, Remark 4.7].)\par
\indent Extending the construction of (2.2), we realise this
solution via a
linear system. We introduce the output space $H_0={\bf
C}^2$, then the Hilbert space $H_1=\ell^2(H_0)$, the state space $H=L^2((t,
\infty );ds; H_1)$ and then let ${\cal D}(A)=\{ f\in H: sf(s)\in H\}$; then we
choose 
$$b_j(s)=\Gamma (jI+W_\infty )^{-1}s^{j-1+W_\infty }\qquad (j=0,
1,\dots ),\eqno(3.8)$$
\noindent recalling that $\Gamma
(z)^{-1}$ is entire. With this choice and some convergence factor $\kappa_0>1$, we introduce linear maps    
$$\matrix{ A:&{\cal D}(A)\rightarrow H:& f(s)\mapsto sf(s);\cr
              B_W:&{}&\beta\mapsto
(\kappa_0^jb_j(s)\beta)_{j=0}^\infty ;\cr
              C:&{\cal D}(A)\rightarrow {\bf C}^2:&
(f_j)_{j=0}^\infty \mapsto \sum_{j=0}^\infty \int_0^\infty
\kappa_0^{-j}C_jf_j(s)\, ds.\cr}\eqno(3.9)$$
\noindent We prove below that $e^{-xA}\beta\in H$ for all sufficiently
large $x$. As usual, we introduce $\Xi_x:L^2(0, \infty )\rightarrow H$ such that
$$\Xi_xf=\int_x^\infty e^{-sA}B_Wf(s)\, ds\eqno(3.10)$$
\noindent and the observability operator $\Theta_x:L^2((0, \infty
);H_0)\rightarrow L^2((t, \infty );H_1)$ by
$$\Theta_x f=\int_x^\infty e^{-sA^\dagger}C_W^\dagger f(s)\,
ds.\eqno(3.11)$$

\vskip.05in
\noindent {\bf Proposition 3.2.} {\sl (i) There exist $\kappa_0,x_0>0$ such
that the operators $\Theta_x:L^2((0,\infty );H_0)\rightarrow H$ and
$\Xi_x: L^2((0, \infty ); H_0)\rightarrow H$ are Hilbert--Schmidt for
$x>x_0$.\par
\indent (ii) For $x>x_0$, the linear system
$(-A,B_W,C_W)$ realises the solution $\Phi$ of (3.1), so that} 
$$\Phi (x;t)=C_We^{-xA}B_W\Phi_0.\eqno(3.12)$$
\indent {\sl (iii) Let $\phi_W (x;t)=C_We^{-xA}B_W$. Then the Hankel
operator on $L^2((x_0, \infty );H_0)$ with symbol $\phi_W $ is trace
class.}
\vskip.05in
\noindent {\bf Proof.} (i) We note that $\Theta_x$ has kernel
$(e^{-su}\kappa_0^{-j}C_j^\dagger )_{j=0}^\infty$, and hence the
Hilbert--Schmidt norm satisfies
$$\eqalignno{\Vert \Theta_x\Vert_{HS}^2&=\sum_{j=0}^\infty
\int_t^\infty \!\!\!\int_x^\infty e^{-2su}\kappa_0^{-2j}\, dsdu\,\Vert
C_j^\dagger\Vert^2_{HS}\cr
&\leq \sum_{j=0}^\infty {{\Vert C_j^\dagger\Vert^2_{HS}e^{-2xt}}\over
{\kappa_0^{2j}4xt}};&(3.13)\cr}$$
\noindent so we choose $\kappa_0$ so that this series converges. For notational convenience, suppose that $\Vert
W_\infty\Vert<1$. Then by the functional equation of $\Gamma$, we have
$$\Vert \Gamma (jI+W_\infty )^{-1}u^{W_\infty +j-1}\Vert \leq
{{u^j\Vert (I+W_\infty )^{-1}\Gamma (I+W_\infty )^{-1}\Vert}\over
{\Gamma (j-1)}}\qquad
(u>1).\eqno(3.14)$$
Next we
observe that $\Xi_x: L^2((x, \infty ); H_0)\rightarrow L^2((t, \infty
); H_1)$ has kernel $(e^{-su}\kappa_0^jb_j(u))_{j=0}^\infty$, and 
hence has Hilbert--Schmidt norm
$$\eqalignno{\Vert \Xi_x\Vert_{HS}^2&=\sum_{j=0}^\infty
\int_x^\infty\!\!\!\int_t^\infty e^{-2su}\kappa_0^{2j}\Vert
b_j(u)\Vert_{HS}^2\, duds\cr
&\leq \sum_{j=0}^\infty \int_t^\infty
\kappa_0^{2j}e^{-2xu}(2u)^{-1}\Vert b_j(u)\Vert^2_{HS}\,
du.&(3.15)\cr}$$
\noindent where the tail of the series is by (3.14)
$$\leq \kappa_W\sum_{j=2}^\infty {{\kappa_0^{2j}\Gamma (2j)}
\over{\Gamma (j-1)^2 (2x)^{2j}}}\Vert (I+W_\infty )^{-1}\Gamma 
(I+W_\infty )^{-1}\Vert^2_{HS} \eqno(3.16)$$
\noindent for some $\kappa_W>0.$ Having chosen $\kappa_0$, we then select
$x_0$ so that the series converges for all $x>x_0;$ then both
$\Theta_x$ and $\Xi_x$ are Hilbert--Schmidt.\par

\indent (ii) Hence we can calculate
$$\eqalignno{C_We^{-xA}B_W&
=\sum_{j=0}^\infty \int_0^\infty C_je^{-xs}b_j(s)\, ds\cr
&=\sum_{j=0}^\infty C_j\Gamma (jI+W_\infty)^{-1}
\int_0^\infty s^{j+W_\infty -1}e^{-sx}\, ds\cr
&=\sum_{j=0}^\infty C_jx^{-W_\infty-j}.&(3.17)\cr}$$  
\indent (iii) By (i), the operator $\Theta_x^\dagger\Xi_x$ is trace
class on $L^2((0, \infty );H_0)$ for all $x>x_0$.\par
\rightline{$\square$}\par
\vskip.05in
\indent Furthermore, the operator
$R_x=\int_x^\infty e^{-sA}B_WC_We^{-sA}\, ds$ on $H$ may be represented as a
kernel with values in a doubly infinite block matrix with $2\times 2$
matrix entries, namely
$$R_x\leftrightarrow \Bigl[
{{\kappa_0^{j-k}b_j(u)C_ke^{-x(u+v)}}\over{u+v}}
\Bigr]_{j,k=0,1, \dots };\eqno(3.18)$$
\noindent this generalises (2.3). Consequently one can in principle
compute the kernel
$$G_W(x,y)=-C_We^{-xA}(I-R_x)^{-1}e^{-yA}B_W,\eqno(3.19)$$
\noindent which satisfies the Gelfand--Levitan equation
$$G_W(x,y)+\phi_W (x+y)-\int_x^\infty G_W(x,w)\phi_W (w+y)\,
dw=0\qquad (t<x<y)\eqno(3.20)$$
\noindent where $\phi_W(x;t)=C_We^{-xA}B_W.$\par

\vskip.05in

\indent We also introduce
$$\sigma_{j,k}=\left[\matrix{I_j&0\cr 0&-I_k\cr}\right]\eqno(3.21)$$
\noindent which has rank $j+k$ and signature $j-k$.\par

\vskip.05in

\noindent {\bf Theorem 3.3.} {\sl Suppose that $W_\infty$ is as in
Lemma 3.1. Let $\Phi (\lambda ;t)$ be a bounded solution of (3.1) 
in $L^2((t,
\infty ); \lambda^{-1}d\lambda ; {\bf R}^2)$ such that $\int_t^\infty
\lambda^{-1}\Vert \Phi (\lambda ;t )\Vert^2d\lambda <\infty$, and let 
$$K(\lambda , \mu ;t)={{\langle J\Phi (\lambda ;t), \Phi (\mu
;t)\rangle }\over{\lambda
-\mu }}.\eqno(3.22)$$
\noindent (i) Then there exists $\phi\in L^2((0,
\infty ); \lambda d\lambda ; {\bf R}^6)$ such that 
$$K(\lambda , \mu ;t)=\int_0^\infty \langle 
\sigma_{3,3}\phi (\lambda +s;t),
\phi (\mu +s;t)\rangle \, ds\qquad (\lambda , \mu >t; \lambda \neq \mu
).\eqno(3.23)$$
\noindent and hence $K$ defines a trace class
operator on $L^2((t, \infty ); d\lambda
).$}\par   
\noindent {\sl (ii) The kernel ${{\partial}\over{\partial t}}K(\lambda
, \mu ;t )$ is of finite rank in $(\lambda , \mu )$.}\par  
\vskip.05in
\noindent {\bf Proof.} Jimbo [14] has shown that the fundamental solution
matrix to (3.1) satisfies
$$Y(x,t)=\bigl( 1+O(x^{-1})\bigr)\left[\matrix {x^{-\theta_\infty/2
}&0\cr 0& x^{\theta_\infty/2 }\cr}\right];\eqno(3.24)$$
\noindent hence there exist solutions that satisfy the hypotheses.\par
\indent (i) We suppress the parameter $t$ to simplify
notation. From the differential equation (3.1), we have
$$\Bigl( {{\partial}\over{\partial \lambda}}+{{\partial }\over{\partial
\mu}}\Bigr) {{\langle J\Phi (\lambda ), \Phi (\mu )\rangle}\over{\lambda
-\mu}}=\Bigl({{1}\over{\lambda -\mu}}\Bigr)\sum_{\nu =0,1,t}
\Bigl\langle \Bigl( {{JW_\nu}\over{\lambda
-\nu}}+{{W_\nu^\dagger J}\over{\mu -\nu}}\Bigr) \Phi (\lambda ), 
 \Phi (\mu )\Bigr\rangle .\eqno(3.25)$$
\indent Now
$$JW_\nu =\left[\matrix{-(z_\nu
+\theta_\nu )/u_\nu& z_\nu+\theta_\nu /2\cr z_\nu +\theta_\nu/2&-u_\nu
z_\nu\cr}\right]\qquad (\nu =0,1,t)\eqno(3.26)$$
\noindent which have rank two and signature zero since $\det
W_\nu =-\theta_\nu^2/4<0$. Hence $JW_\nu =V_\nu^\dagger
\sigma_{1,1}V_\nu$ for some $2\times 2$ real matrix $V_\nu ,$ and
$JW_\nu=V_\nu^\dagger\sigma_{1,1}V_{\nu}$. Thus we find that (3.25) reduces to
$$-{{\langle \sigma_{1,1}V_0\Phi (\lambda ), V_0\Phi (\mu
)\rangle}\over{\lambda\mu}} -{{\langle \sigma_{1,1}V_1\Phi (\lambda ),
V_1\Phi (\mu
)\rangle}\over{(\lambda -1)(\mu -1)}} -{{\langle \sigma_{1,1}
V_t\Phi (\lambda ), V_t\Phi (\mu
)\rangle}\over{(\lambda-t)(\mu-t)}}.\eqno(3.27)$$
\noindent Let  
$$\phi (\lambda )=\left[ \matrix{ {{V_0\Phi (\lambda
)}\over{\lambda}}\cr 
{{V_1\Phi (\lambda )}\over{\lambda -1}}\cr
{{V_t\Phi (\lambda )}\over{\lambda -t}}\cr}\right],\eqno(3.28)$$
\noindent which satisfies, after we permute the coordinates in the
obvious way, 
$$\eqalignno{-\sum_{\nu =0,1,t} {{\langle \sigma_{1,1}
V_\nu \Phi (\lambda ), V_\nu \Phi (\mu
)\rangle}\over{(\lambda -\nu )(\mu -\nu )}}&=-\langle \sigma_{3,3}\phi (\lambda ),
\phi (\mu )\rangle\cr
&=\Bigl( {{\partial }\over{\partial \lambda}}+{{\partial
}\over{\partial\mu}}\Bigr) \int_0^\infty \langle \sigma_{3,3}\phi (\lambda +s ), \phi
(\mu +s)\rangle \, ds.&(3.29)\cr}$$
\noindent We observe that both sides of (3.23) converge to zero as
$\lambda\rightarrow\infty$ and as $\mu\rightarrow\infty$. 
By comparing the derivatives as in (3.25) and (3.29), we deduce
(3.23).\par  
\indent Then $K=\Gamma_\phi^\dagger\sigma_{3,3}\Gamma_\phi .$
\noindent We observe that the Hilbert--Schmidt norm of $\Gamma_\phi$
satisfies
$$\eqalignno{\Vert\Gamma_\phi\Vert_{HS}^2&=\int_{t}^\infty (\lambda
-t)\Vert\phi (\lambda )\Vert^2\, d\lambda\cr
&\leq \kappa \int_{t}^\infty {{\Vert \Phi (\lambda
)\Vert^2}\over{\lambda}}\, d\lambda&(3.30)\cr}$$
\noindent for some $\kappa >0$, so $K$ gives a trace class operator on $L^2(t, \infty
)$.\par  
\indent (ii) By a similar calculation, one can compute the derivative of $K$
with respect to the position of the critical point, and find
$$\eqalignno{{{\partial }\over{\partial t}}K(\lambda ,\mu
;t)&={{1}\over{(\lambda -t)(\mu -t)}}\Bigl\langle \left[\matrix{
-(z_t+\theta_t)/u_t& z_t+\theta_t/2\cr z_t+\theta_t/2&-u_tz_t\cr}\right]
\Phi (\lambda ;t), \Phi (\mu ;t)\Bigr\rangle ;&(3.31)\cr}$$
\noindent evidently this is a finite sum of products of functions of
$\lambda$ and functions of $\mu$ for each $t$.\par
\rightline{$\square$}\par

\vskip.05in
\noindent {\bf 4. The $\tau$ function associated with Painlev\'e's
equation VI}\par 
\vskip.05in 
\indent In [2], Ablowitz and Segur derived an integral equation  involving the Airy kernel for
the solutions of $P_{II}$. Here we solve an integral equation and
derive an expression for $\det (I-KP_{(x, \infty )})$, which is associated with $P_{VI}$. From Proposition 3.2, we recall the linear system $(-A_W, B_W, C_W)$
that realises $\phi_W$, and likewise we introduce a linear system
$(-A_V, B_V,C_V)$ that realises $\phi_V={\hbox{diagonal}}(V_0/x,
V_1/(x-1), V_t/(x-t))$;  then by considering 
$$(-(A_V\otimes I+I\otimes A_W), B_V\otimes B_W, C_V\otimes C_W)$$
\noindent we introduce a new linear system that realises $\phi$ from
Theorem 3.3, so that $\phi (x)=Ce^{-xA}B.$\par
\indent Next we let $\Gamma_\phi$ be the Hankel
integral operator with symbol $\phi $; also let
$\phi_{(x)}(y)=\phi (y+2x)$ and let $L_x$ be observability Gramian
$$L_x=\int_x^\infty e^{-sA}BB^\dagger e^{-sA^\dagger}\,
ds=\Xi_x\Xi_x^\dagger .\eqno(4.1)$$
\noindent To take account of the signature, we introduce the
 the modified controllability Gramian
$$Q^\sigma_x=\int_x^\infty e^{-sA^\dagger}C^\dagger
\sigma_{3,3}Ce^{-sA}\, ds.\eqno(4.2)$$
\indent We also introduce the $(6+1)\times (6+1)$ block matrices
$$ G(x,y)=\left[\matrix{U(x,y)&V(x,y)\cr T(x,y)&\zeta
(x,y)\cr}\right]\eqno(4.3)$$

\noindent and
$$\Phi (x)=\left[\matrix{ 0&\phi(x)\cr \phi (x)^\dagger
&0\cr}\right],\eqno(4.4)$$
\noindent and the Gelfand--Levitan integral equation 
$$G(x,y)+\Phi (x+y)+\int_x^\infty G(x,w)\ast \Phi (w+y)\,
dw=0,\eqno(4.5)$$
 \noindent where we have introduced a special matrix product to
incorporate the signature, namely 
$$\eqalignno{&\int_x^\infty G(x,w)\ast\Phi
(w+y)\,dw\cr
&\quad =\left[\matrix{\int_x^\infty V(x,w)\phi (w+y)^\dagger \sigma_{3,3}
dw& \int_x^\infty U(x,w)\phi (w+y)dw\cr 
\int_x^\infty \zeta (x,y)\phi (w+y)^\dagger dw&
\int_x^\infty T(x,w)\sigma_{3,3}\phi (w+y)dw\cr}\right].&
(4.6)\cr}$$

\vskip.05in
\noindent {\bf Theorem 4.1.} {\sl Suppose that $Q_x$ and
$L_x$ are trace-class operators with operator norms less than one
for all $x>t$.
Then there exists a solution to the integral equation (4.5)
such that $\tau_K(x)=\det (I-P_{(x,\infty )}K)$ satisfies}
$${{d}\over{dx}}\log
\tau_K(x)={\hbox{trace}}\,G(x,x).\eqno(4.7)$$
\vskip.05in
\noindent {\bf Proof.} By Theorem 3.3, we have
$K=\Gamma_\phi^\dagger\sigma_{3,3}\Gamma_\phi $, and so 
$$\eqalignno{ \tau_K(x)&=\det (I-P_{(x, \infty
)}\Gamma_\phi^\dagger\sigma_{3,3}\Gamma_\phi)\cr
&=\det (I-\Xi_x^\dagger \Theta_x\sigma_{3,3}\Theta_x^\dagger\Xi_x)\cr
&=\det (I-
\Theta_x\sigma_{3,3}\Theta_x^\dagger \Xi_x\Xi_x^\dagger)\cr
&=\det (I-Q^\sigma_xL_x).&(4.8)\cr}$$
\noindent One can verify that  
$$\eqalignno{{}&\left[\matrix{U(x,y)&V(x,y)\cr T(x,y)&\zeta
(x,y)\cr}\right]&(4.9)\cr
&\qquad  =\left[ \matrix{ Ce^{-xA}(I-L_xQ_x^\sigma
)^{-1}L_xe^{-yA^\dagger} C^\dagger \sigma_{3,3}  &
-Ce^{-xA}(I-L_xQ_x^\sigma)^{-1}e^{-yA}B\cr
-B^\dagger e^{-xA^\dagger}(I-Q_x^\sigma
L_x)^{-1}e^{-yA^\dagger}C^\dagger &B^\dagger e^{-xA^\dagger}
(I-Q_x^\sigma L_x)^{-1}Q_x^\sigma
e^{-yA}B\cr}\right]}$$  
\noindent gives a solution to (4.6), so that
$$\eqalignno{{\hbox{trace}}\,U(x,x)&={\hbox{trace}}\Bigl(
(I-L_xQ^\sigma_x)^{-1}L_x e^{-xA^\dagger}C^\dagger
\sigma_{3,3}Ce^{-xA}\Bigr)\cr
&=-{\hbox{trace}}\Bigl(
(I-L_xQ^\sigma_x)^{-1}L_x{{dQ^\sigma_x}\over{dx}}\Bigr).&(4.10)\cr}$$
\noindent Likewise we have
$$\eqalignno{\zeta (x,x)&={\hbox{trace}}\Bigl(
(I-Q_xL_x)^{-1}Q^\sigma_x e^{-xA}BB^\dagger e^{-xA^\dagger}\Bigr)\cr
&=-{\hbox{trace}}\Bigl(
(I-Q^\sigma_x
L_x)^{-1}Q^\sigma_x{{dL_x}\over{dx}}\Bigr).&(4.11)\cr}$$
\indent Adding and rearranging, we obtain
$$\eqalignno{ {\hbox{trace}}\,G(x,x)&=
\zeta (x,x)+{\hbox{trace}}\, U(x,x)\cr
&=-{\hbox{trace}}\Bigl(
(I-L_xQ^\sigma_x)^{-1}L_x{{dQ^\sigma_x}\over{dx}}\Bigr)\cr
&\qquad -{\hbox{trace}}\Bigl(
(I-L_xQ^\sigma_x)^{-1}{{dL_x}\over{dx}}Q^\sigma_x\Bigr)\cr
&={{d}\over{dx}}{\hbox{trace}} \log (I-L_xQ^\sigma_x)\cr
&={{d}\over{dx}}\log \tau_K(x).&(4.12)\cr}$$  
\rightline{$\square$}\par
\vskip.05in

\indent We introduce
the new variable $u$ by the elliptic integral
$$u(y,t)=\int_\infty^y {{d\lambda}\over{\sqrt{\lambda (\lambda -1)(\lambda
-t)}}},\eqno(4.13)$$
\noindent then we let $Z={{\partial y}\over{\partial u}}$ and
$Y=y$, so $(Y,Z)$ lies on the elliptic curve $Z^2=Y(Y-1)(Y-t)$ which
depends upon the parameter $t$. Soon after his discovery of $P_{VI}$, 
 R. Fuchs showed that if $y(t)$ satisfies $P_{VI}$, then
$u(t)=u(y(t),t)$ 
satisfies  
$$\eqalignno{-t(1-t)&{{d^2u}\over{dt^2}}+(2t-1){{du}\over{dt}}+{{u}\over{4}}
\cr
&=-{{\sqrt {y(y-1)(y-t)}}\over{t(1-t)}}\Bigl( 2\alpha +{{2\beta
t}\over{y^2}}+{{\gamma (t-1)}\over{(y-1)^2}}+(\delta
-1/2){{t(t-1)}\over{(y-t)^2}}\Bigr),&(4.14)\cr}$$
\noindent where we recognise Legendre's differential operator on the
left-hand side; see [32, p 304]. By analysing these solutions, Guzzetti [13]
obtains various series representations and bounds on the growth of
$y(t)$. We can analyse symbols that are elliptic
functions of the second kind since their transfer functions have
special properties.\par

\vskip.1in

\noindent {\bf 5. Kernels associated with the hypergeometric
equation}\par
The $P_{VI}$ equation is closely related to Gauss's
hypergeometric equation [32, p 283]
$$\lambda (1-\lambda ){{d^2f}\over{d\lambda^2}}+(c-(a+b+1)\lambda )
{{df}\over{d\lambda }}-abf(\lambda )=0.\eqno(5.1)$$
\noindent We introduce $c_0=c$ and $c_1=a+b-c+1$, then introduce the matrix
$$W(\lambda )=\left[\matrix{0& \lambda^{-c_0}(\lambda -1)^{-c_1}\cr
-ab\lambda^{c_0-1}(\lambda -1)^{c_1-1}&0\cr}\right]\eqno(5.2)$$
\noindent so that we can express (5.1) in the form of a first
order linear differential equation as in (5.4). For special choices
of the parameters $a,b,c$, we can obtain a factorization of the
corresponding kernel (5.5) which has the form of (1.2). For a separable Hilbert space $H$ we introduce the identity
operator $I_H$ and 
$$\sigma_{H,H}=\left[\matrix{I_H&0\cr 0&
-I_H\cr}\right].\eqno(5.3)$$
\vskip.05in
\noindent {\bf Theorem 5.1.} {\sl Suppose that $0\leq c\leq 1$ and
$a+b=0$, that $2\sqrt{-ab}$ is not an integer, and that $-ab>5/4,$ and let $\Psi$ be a bounded solution for the equation
$${{d\Psi}\over{d\lambda }}=W(\lambda )\Psi (\lambda ),\eqno(5.4)$$
\noindent such that $\int_1^\infty x\Vert\Psi (x)\Vert^2dx<\infty$; then let 
$$K(x,y)={{\langle J\Psi (x), \Psi (y)\rangle}\over{x-y}}\qquad (x\neq
y; x,y>1).\eqno(5.5)$$
\noindent (i) Then there exists a separable Hilbert space $H$ and 
$\phi :(1, \infty )\rightarrow H^2$
such that\par
\noindent  $\int_{1+\delta }^\infty
x\Vert \phi
(x)\Vert^2_{H^2}\,
dx<\infty$ and $K=\Gamma_\phi^\dagger\sigma_{H,H}\Gamma_\phi$ so that $K$
defines a trace class kernel on $L^2((1+\delta , \infty ); dx)$
for all $\delta >0$.}\par
\noindent {\sl 
\noindent (ii) The statement of Theorem 4.1 applies to 
$$\tau_K(s)=\det (I-KP_{(s, \infty )})=\det
(I-\Gamma_{\phi_{(s/2)}}^\dagger\sigma_{H,H}\Gamma_{\phi_{(s/2)}}),\eqno
(5.6)$$
\noindent with
obvious changes to notation; so ${{d}\over {dt}}\log\tau_K (t)$ is
given by the diagonal of the solution of a Gelfand--Levitan
equation.\par
\noindent (iii) If moreover $c$ is rational, then $K$ arises from a
Tracy--Widom system as in (1.1).\par

}\par 
\vskip.05in
\noindent {\bf Proof.} Let 
$$q(\lambda )={{-ab}\over{\lambda (\lambda
-1)}}+{{1}\over{4}}\Bigl({{c^2-2c}\over{\lambda^2
}}+{{2c(1-c)}\over{\lambda (\lambda -1)}}+{{c^2-1}\over{(\lambda
-1)^2}}\Bigr),\eqno(5.7)$$
\noindent which is asymptotic to $(-ab-1/4)/\lambda^2$ as
$\lambda\rightarrow\infty$. By the
Liouville--Green transformation [25,
p.229] , we
can obtain solutions to (5.1) with asymptotics of the form
$$f_{\pm}(\lambda )\asymp \lambda^{-c/2}(\lambda -1)^{-(1-c)/2}q(\lambda
)^{-1/4}\exp\Bigl( \pm \int_2^\lambda q(x)^{1/2}\,
dx\Bigr)\qquad (\lambda \rightarrow\infty ),\eqno(5.8)$$
\noindent and one can deduce that 
$\int_2^\infty xf_-(x)^2\, dx <\infty $. Hence there exist solutions
that satisfy the hypotheses.\par
\indent (i) 
We observe that $c_1+c_0=1$,
so $0\leq c_0, c_1, 1-c_0, 1-c_1\leq 1$; we assume that
$0<c_0,c_1<1$, as the cases of equality are easier. Evidently
the functions $\lambda^{-c_0}(\lambda
-1)^{c_0-1}$ and $\lambda^{-c_1}(\lambda -1)^{c_-11}$ are 
operator monotone
decreasing on $(1,\infty )$ in Loewner's sense and by [1, p.577] we
have an integral representation
$$\lambda^{-c_0}(\lambda -1)^{c_0-1}={{\sin \pi
c_0}\over{\pi}}\int_{-1}^0{{(-u)^{-c_0}(1+u)^{c_0-1}\, du}\over{
\lambda +u}}\qquad (\lambda >1);\eqno(5.9)$$
 clearly a similar representation holds for 
$\lambda^{-c_1}(\lambda -1)^{c_1-1}$ with $c_1$ instead of $c_0$. Hence there exist positive
measures $\omega_1 $ and $\omega_0$ on $[-1,0]$ such that 
$${{JW(x)+W(y)^\dagger
J}\over{x-y}}=\left[\matrix{ab{{x^{-c_1}(x-1)^{c_1-1}
-y^{-c_1}(y-1)^{c_1-1}}\over{x-y}}&0\cr 
0&{{x^{-c_0}(x-1)^{c_0-1}-y^{-c_0}(y-1)^{c_0-1}}\over{x-y}}\cr}\right]$$
$$=\int_{-1}^{0}{{1}\over{(x+u)(y+u)}}\left[\matrix{-ab\omega_1
(du)&0\cr 0& -\omega_0 (du)\cr}\right] \eqno(5.10)$$  
\noindent in which $-ab\geq 0$. The matrix kernel $(JW(x)+W(y)^\dagger
J)/(x-y)$ operates as a Schur multiplier on the rank one tensor $\Psi
(x)\otimes \Psi (y)$ in $L^2((1+\delta , \infty ); {\bf R}^2)$; hence 
for each $\delta >0$, there exists $\kappa_\delta >0$
such the Schur multiplier norm is bounded by $\kappa_\delta $. Since
$\Psi (x+s)$ gives a Hilbert--Schmidt kernel, the 
operator $\int_0^\infty \Psi (x+s)\otimes \Psi (y+s)\, ds$ is
trace class on $L^2((1+\delta , \infty ); dx)$, and it follows that 
$$K(x,y)=\int_0^\infty \Bigl\langle {{JW(x+s)+W(y+s)^\dagger
J}\over{x-y}}\Psi (x+s), \Psi (y+s)\Bigr\rangle ds \eqno(5.11)$$
\noindent is also trace class. 
As in  
Theorem 1.1 of [4], we can introduce the Hilbert space $H$, $\phi\in
L^2((1+\delta ,\infty ); xdx; H^2)$ and the
 Hankel operator $\Gamma_\phi$ with symbol $\phi$ such that 
$K=\Gamma_\phi^\dagger
\sigma_{H,H}\Gamma_\phi ,$ so 
$$K(x,y)=\int_0^\infty \langle  \sigma_{H,H}\phi (x+s), \phi
(y+s)\rangle_{H^2}\, ds\eqno(5.12)$$
where $\sigma_{H,H}$ takes account of the fact that
the Schur multiplier is positive on the top left matrix block and 
negative on the bottom right matrix block.\par
\indent (ii) We observe that 
$$W(\lambda )={{1}\over{\lambda}}\left[ \matrix{ 0&1\cr
-ab&0\cr}\right] +O(\lambda^{-2})\qquad
(\vert\lambda\vert\rightarrow\infty ),\eqno(5.13)$$
\noindent is analytic at infinity and the residue matrix has
eigenvalues $\pm \sqrt {-ab}$ which do not differ by a positive integer.
Hence we can repeat the proof of Lemma 3.1 and realise
the solution $\Psi $ of (5.4) by a linear system involving the
coefficients in the Laurent series of $\Psi$. Then we can realise
$\phi \in L^2((0, \infty ); H^2)$ by means of a linear system $(-A, B,C)$, where the state space
is $L^2((0, \infty ); H^2)$. We can now follow through the proof in
section 4 as express $\tau$ in terms of the Gelfand--Levitan
equation.\par 
\indent (iii) Let $c=k/n$; then $\{ (X,Z): Z^n=X^k(X-1)^{n-k}\}$ gives a
$n$-sheeted cover of ${\bf P}^1$, ramified at $0,1,\infty$. On this
compact Riemann
surface, the functions $\lambda^{-c_0}(\lambda -1)^{c_0-1}$ and
$\lambda^{-c_1}(\lambda -1)^{c_1-1}$ are rational.\par

\rightline{$\square$}\par  
\vskip.05in
\noindent {\bf Remarks.} (i) The Painlev\'e equations can be expressed as Hamiltonian
systems in the canonical variables $(\lambda ,\mu )$, where the
Hamiltonian is a rational function of $(\lambda , \mu )$; see [24] for
a list. Okamoto [24] showed that there exists a holomorphic
function $\tau$ on the universal covering surface of ${\bf P}^1\setminus
\{ 0,1,\infty\}$ such that $H_{VI}(t, \lambda (t), \mu (t))={{d}\over {dt}}\log \tau
(t).$ The methods of [11, 
15, 16] involve complex analysis and
differential geometry, and are not intended to address the
operator properties of $K$.\par
\indent (ii) Borodin and Deift [7] have identified an integrable kernel $K$
involving solutions ${}_2F_1$ of the hypergeometric equation and
considered 
$\tau (t)=\det (I-P_{(t, \infty )}K)$; they showed that $\sigma
(t)={{d}\over{dt}}\log \tau (t)$ satisfies the Jimbo--Miwa
$\sigma$ form of $P_{VI}$.\par


\vskip.05in
\noindent {\bf 6. The $\tau$ function associated with a Hankel
operator on exponential bases}\par
\vskip.05in
\noindent We
wish to find a more explicit expression for $\tau $ and for 
$\sigma (t)={{d}\over{dt}}\log \tau (t)$ for suitable $K$, especially those
$K$ that factor as $K=\Gamma_\phi^\dagger\Gamma_\phi$. We can obtain an explicit formula for $\tau$ when $\phi$ has the
exponential expansion 
$$\phi (x)=\sum_{j=1}^\infty \xi_j e^{-\lambda_jx}\eqno(6.1)$$
\noindent where the coefficients $\xi_j$ lie in some Hilbert space $E$.
In this section we establish the existence of such expansions by using
the theory of approximation of compact Hankel operators, whereas in
subsequent sections we consider the transfer function $\hat \phi (s)$ of 
$\phi$ and use the Mittag-Leffler expansion to give explicit
formulas. The Hankel operator with symbol $\phi$ can be expressed in
terms of the exponential basis as a
relatively simple matrix, so we can derive expressions for its Fredholm
determinant. Our
applications in sections 7 and 8 are to cases in which the poles lie on an arithmetic
progression, which occurs when $\phi$ is a theta function or arises
by a certain transformation of a power series.\par
\indent We suppose that $\lambda_j\in {\bf C}$ with $\Re
\lambda_j>0$ are such that
$(e^{-t\lambda_j})_{j=1}^\infty$ are linearly independent exponentials,
so that
$$D_N=\det\Bigl[
{{1}\over{\lambda_j+\bar\lambda_k}}\Bigr]_{j,k=1}^N>0\qquad (N=1, 2,
\dots ).\eqno(6.2)$$
\indent Suppose that $\xi =(\xi_j)_{j=1}^\infty \in\ell^1$ and
introduce the operators
$$\matrix{B:&{\bf C}\rightarrow\ell^1\subset \ell^2:
 &a\mapsto a\xi\cr 
              e^{-tA}:&\ell^2\rightarrow\ell^2: &(\alpha_j)_{j=1}^\infty
\mapsto (e^{-t\lambda_j}\alpha_j)_{j=1}^\infty\cr
C:&\ell^1\subset\ell^2\rightarrow {\bf C}:&
(\alpha_j)_{j=1}^\infty\mapsto \sum_{j=1}^\infty
\alpha_j\cr
\Theta :&L^2(0, \infty )\rightarrow \ell^2: & f\mapsto
(\int_0^\infty e^{-\bar\lambda_j s}f(s)\,
ds)_{j=1}^\infty.\cr}\eqno(6.3)$$
\vskip.05in

\noindent {\bf Theorem 6.1.} {\sl Suppose that $\Theta$ is bounded
and that there exist constants
$\delta , M>0$ such that $\Re\lambda_j\geq \delta$ and
$\sum_{k=1}^\infty \vert \lambda_j+\lambda_k\vert^{-2}\leq M$ for
all $j$; let $\xi\in \ell^1$.\par
\indent (i) Then the symbol $\phi (x)=Ce^{-xA}B$ gives rise to a Hankel operator
$\Gamma_\phi:L^2(0, \infty )\rightarrow L^2(0, \infty )$ which is trace class.\par
\indent (ii) The operator 
$$R_x=\int_x^\infty e^{-sA}BCe^{-sA}\, ds\eqno(6.4)$$
\noindent on $\ell^2$ is trace class, and for $\mu$ is an open
neighbourhood of zero, the kernel $T_\mu
(x,y)=-\mu Ce^{-xA}(I+\mu R_x)^{-1}e^{-yA}B$ gives a solution to the integral
equation
$$T_\mu (x,y)+\mu \phi (x+y) +\mu\int_x^\infty T_\mu (x,z)\phi (z+y)\,
dz=0\qquad (0<x\leq y).\eqno(6.5)$$
\indent (iii) Suppose that $(I-R_t)$ is invertible for all $t>0$.
Then the Hankel operator $\Gamma_{\phi_{(t)}}$ with kernel 
$\phi (x+y+2t)$ satisfies}
$$\det (I-\Gamma_{\phi_{(t)}})=\exp\Bigl( -\int_t^\infty T_{-1}(u,u)\,
du\Bigr).\eqno(6.6)$$
\vskip.05in
\noindent {\bf Proof.} (i) The kernel may be expressed as a sum of rank-one kernels
$$\Gamma_\phi \leftrightarrow\sum_{j=1}^\infty
\xi_je^{-\lambda_j(x+y)}\eqno(6.7)$$
\noindent where $\sum_{j=1}^\infty \vert\xi_j\vert/\Re\lambda_j$
converges, so $\Gamma_\phi$ is trace class.\par
\indent (ii) By considering the rows of the matrix
$$R_x\leftrightarrow
\Bigl[{{\xi_je^{-(\lambda_j+\lambda_k)x}}\over{\lambda_j+\lambda_k}}
\Bigr]_{j,k=1}^\infty\eqno(6.8)$$
\noindent we see that $R_x$ is also trace class. When $\vert\mu\vert\Vert R_x\Vert <1$, the kernel
$T_\mu (x,y)$ is well defined, and one verifies the identity (6.5) by
substituting.\par
\indent (iii) The operators 
$$C:\ell^1\rightarrow {\bf C}, \qquad e^{-tA}:\ell^1\rightarrow\ell^1,
\qquad R_x:\ell^1\rightarrow \ell^1,\qquad B:{\bf C}\rightarrow
\ell^1\eqno(6.9)$$
\noindent are all bounded, and $\xi\mapsto R_x$ is continuous from
$\ell^1$ to the trace class; hence $T(x,y)$ depends continuously on
$\xi$ in a neighbourhood of $0$ in $\ell^1$. Suppose that 
$(\xi^{(n)})_{n=1}^\infty$
is a sequence of vectors in $\ell^1$ that have only finitely many
nonzero terms, and that $\xi^{(n)}\rightarrow\xi$ as
$n\rightarrow\infty$. Denoting the operators corresponding to
$\xi^{(n)}$ by $R_x^{(n)}$ etcetera, we can manipulate the finite
matrices and deduce that
$$T_{-1}^{(n)}(x,x)={{d}\over{dx}}\log\det (I-R_x^{(n)})\eqno(6.10)$$
\noindent and hence
$$\int_s^t T_{-1}^{(n)}(x,x)\, dx=\log\det (I-R_t^{(n)}) -
\log\det (I-R_s^{(n)});\eqno(6.11)$$
\noindent so letting $n\rightarrow\infty$, we deduce that
 $$\int_s^t T_{-1}(x,x)\, dx=\log\det (I-R_t) -
\log\det (I-R_s).\eqno(6.12)$$
\indent The operator $\Xi :L^2(0, \infty )\rightarrow \ell^2$ 
given by
$$\Xi f=\int_0^\infty e^{-tA}Bf(t)\, dt\eqno(6.13)$$
\noindent has matrix representation
$$\Xi\Xi^\dagger \leftrightarrow\Bigl[
{{\xi_j\bar\xi_k}\over{\lambda_j+\bar \lambda_k}}\Bigr]_{j
,k=1}^\infty\eqno(6.14)$$
\noindent with respect to the standard basis $(e_j)$, and hence $\Xi$ is
Hilbert--Schmidt since \par
\noindent $\sum_{j=1}^\infty \Vert\Xi^\dagger  e_j\Vert^2
<\infty $. The operator $\Theta$ is bounded by hypothesis, hence
$\Theta^\dagger$ is also bounded;  so
$R_0=\Xi\Theta^\dagger$ is also Hilbert--Schmidt.\par  

\indent The operator $\Gamma_\phi$ is trace 
class by (ii), and 
the non-zero eigenvalues of $\Gamma_\phi =\Theta^\dagger\Xi$ and
$R_0=\Xi\Theta^\dagger$ are equal, hence
$$\det (I-\Gamma_{\phi_{(x)}})=\det (I-R_x)\eqno(6.15)$$
\noindent which when combined with (6.12), implies that 
$$\log\det (I-
\Gamma_{\phi_{(s)}})-\log\det (I-
\Gamma_{\phi_{(t)}})=\int_t^sT_{-1}(u,u)\, du.\eqno(6.16)$$
\noindent Evidently $\Gamma_{\phi_{(s)}}\rightarrow 0$ as
$s\rightarrow\infty$, and hence (6.6) follows from (6.16).\par
\rightline{$\square$}\par
\vskip.1in

\vskip.05in

\noindent {\bf Theorem 6.2.} {\sl Let $K$ be an integral operator on
$L^2((0,\infty ); dt; {\bf C})$ such that:\par
\indent (i) $0\leq K\leq I$ and $I-K$ is invertible;\par
\indent (ii) there exists a separable Hilbert space $E$ and $\phi\in
L^2((0, \infty ); tdt;E)$ such that
$K=\Gamma_\phi^\dagger\Gamma_\phi$.\par
\noindent Then $K$ has a $\tau$-function $\tau_K$ and there exists a 
sequence $(K_n)_{n=1}^\infty$ of finite rank integral
operators with corresponding $\tau$-functions $\tau_{K_n}$ such
that:\par
\indent (1) $K_n\rightarrow K$ in trace class norm;\par
\indent (2) $\tau_{K_n}(x)\rightarrow \tau_K(x)$ uniformly on compact
sets as $n\rightarrow\infty;$\par
\indent (3) $\tau_{K_n}(x)=\sum_{j=1}^{N_n}a_{jn}e^{-\mu_{jn}x}$ for
some $a_{jn}, \mu_{jn}\in {\bf C}$ with $\Re \mu_{jn}>0$ that are given
in Proposition 6.4 below.}\par
\vskip.05in
\noindent {\bf Proof.} (1) For $\phi\in L^2((0, \infty );tdt;
E)$, the operator $\Gamma_\phi$ is Hilbert--Schmidt and hence $K$ is 
trace class. By the Adamyan--Arov--Krein theorem [26], there exists a
sequence $(\Gamma_{\phi^{(n)}})_{n=1}^\infty$ of finite-rank Hankel
operators such that $\Gamma_{\phi^{(n)}}\rightarrow \Gamma_\phi$ in
Hilbert--Schmidt norm.\par
\indent  Kronecker showed that a Hankel
operator $\Gamma_{\phi^{(n)}}$ has finite rank if and only if the transfer 
function $\hat \phi^{(n)} (s)$ is
rational; see [26]. Hence the typical form for $\phi^{(n)}$ is a finite sum
$$\phi^{(n)}(t)=\sum_{j,k} \xi_{k,j}t^ke^{-\lambda_jt}\eqno(6.17)$$
\noindent where $\xi_{k,j}\in E$ and $\Re \lambda_j>0$; the terms with 
factor $t^k$ give poles of order $k+1$. To resolve the poles of order
greater than one into sums of simple poles, we introduce the difference operator
$\Delta_\varepsilon$ by $\Delta_\varepsilon g(\lambda
)=\varepsilon^{-1}(g(\lambda +\varepsilon )-g(\lambda ))$, which
satisfies $\lim_{\varepsilon\rightarrow 0}\Delta_\varepsilon^kg(\lambda )=g^{(k)}(\lambda )$
whenever $g$ is $k$-times differentiable with respect to $\lambda $. By
the dominated convergence theorem, 
$$\int_0^\infty t\vert k!(-\Delta_\varepsilon )^k
e^{-\lambda_jt}-t^ke^{-\lambda_jt}\vert^2\, dt\rightarrow 0\eqno(6.18)$$
\noindent as $\varepsilon \rightarrow 0$, so we can replace $t^ke^{-\lambda_jt}$
by $k!(-\Delta_\varepsilon)^k e^{-\lambda_jt}$ at the cost of a small
change in the operator $\Gamma_{\phi^{(n)}}$ in Hilbert--Schmidt norm.
Thus we eliminate poles of order greater than one, and we can
ensure that $0\leq
\Gamma_{\phi^{(n)}}^\dagger\Gamma_{\phi^{(n)}}\leq I$, with
$I-\Gamma_{\phi^{(n)}}^\dagger\Gamma_{\phi^{(n)}}$ invertible. Let
$K_n=\Gamma_{\phi^{(n)}}^\dagger\Gamma_{\phi^{(n)}}$ so that $K_n$
has finite rank and $K_n\rightarrow K$ as in trace norm as $n\rightarrow \infty$.\par

\indent (2) Let
$\phi_{(x)}(t)=\phi (t+2x)$ and
$\phi_{(x)}^{(n)}(t)=\phi^{(n)}(t+2x)$. We have
$\Gamma_{\phi^{(n)}_{(x)}}^\dagger\Gamma_{\phi^{(n)}_{(x)}}\rightarrow
\Gamma_{\phi_{(x)}}^\dagger \Gamma_{\phi_{(x)}}$ in trace class norm as
$n\rightarrow\infty$ so  
$$\eqalignno{\tau (x)&=\det (I-KP_{(x, \infty )})\cr
&=\det (I-\Gamma_{\phi_{(x)}}^\dagger\Gamma_{\phi_{(x)}})\cr
&=\lim_{n\rightarrow\infty }\det (I-\Gamma_{\phi_{(x)}^{(n)}}^\dagger
\Gamma_{\phi_{(x)}^{(n)}})\cr
&=\lim_{n\rightarrow\infty}\tau_{K_n}(x)&(6.19)}$$
since the Fredholm determinant is a continuous
functional on the trace class operators.\par 
\indent (3) To calculate the function $\tau_{K_n}(x)$ in (3) of 
Theorem 6.2, we assume that $\phi^{(n)}$ has the form
$$\phi^{(n)} (t)=\sum_{j=1}^N \xi_j^\dagger e^{-\bar\lambda_jt}\qquad
(t>0)\eqno(6.20)$$ 
\noindent where $\xi_j\in E$ and $\Re \lambda_j>0$.  Without loss of
generality we can replace $E$ by the subspace ${\hbox{span}} (
\xi_j)_{j=1}^N$ and for notational simplicity we take $\xi_j\in
M_{1,\nu }({\bf C})$ where $\nu\leq N.$
\vskip.05in
\noindent We introduce
$$a_j={\hbox{row}}\Bigl[ {{\xi_j e^{-2\lambda_jx}}\over
{\lambda_j+\bar\lambda_k}}\Bigr]\in M_{1, \nu N}({\bf C})\eqno(6.21)$$
\noindent and
$$b_m={\hbox{column}}\Bigl[{{\xi_k^\dagger
e^{-2\bar\lambda_kx}}\over{\bar\lambda_k+\lambda_m}}\Bigr]_{k=1}^N\in
M_{\nu N,1}({\bf C}).\eqno(6.22)$$
\rightline{$\square$}\par
\vskip.05in
\noindent 
\noindent {\bf Lemma 6.3.} {\sl The matrix
$$K=[a_jb_m]_{j,m=1}^N\eqno(6.23)$$
\noindent represents the operator
$\Gamma^\dagger_{\phi^{(n)}_{(x)}}\Gamma_{\phi^{(n)}_{(x)}}$ with respect to the
(non-orthogonal) basis $(e^{-\lambda_js})_{j=1}^N.$}\par
\vskip.05in
\noindent {\bf Proof.} We observe that the transfer function of
$\phi^{(n)}_{(x)}$ is the rational function
$$\hat \phi^{(n)}_{(x)} (s)=\sum_{j=1}^\nu {{\xi_j^\dagger
e^{-2\lambda_jx}}\over{s+\lambda_j}}.\eqno(6.24)$$
\noindent The operator $\Gamma^\dagger_{\phi^{(n)}_{(x)}}\Gamma_{\phi^{(n)}_{(x)}}$ 
has kernel in the variables $(s,t)$ 
$$\int_0^\infty \langle \phi^{(n)} (2x+s+u), \phi^{(n)} (2x+t+u)\rangle \, du
\eqno(6.25)$$
\noindent and hence one computes 
$$\Gamma^\dagger_{\phi^{(n)}_{(x)}}\Gamma_{\phi^{(n)}_{(x)}}:e^{-\lambda_m
s}\mapsto \sum_{j,k=1}^N{{\langle \xi_j,
\xi_m\rangle
e^{-2(\bar\lambda_k+\lambda_j)x}}\over{(\lambda_j+\bar\lambda
_k)(\bar\lambda_k+\lambda_m)}}e^{-\lambda_js}.\eqno(6.26)$$
\noindent Recalling the definitions (6.21) and (6.22), one
computes
$$a_jb_m=\sum_{j=1}^N{{\langle\xi_j,\xi_k\rangle
e^{-2(\lambda_j+\bar\lambda_k) x}}\over
{(\lambda_j+\bar\lambda_k)(\bar\lambda_k+\lambda_m)}}\eqno(6.27)$$
\noindent and by comparing this with (6.23), one obtains the stated
identity.\par
\rightline{$\square$}\par
\vskip.05in

\noindent We can proceed to compute the $\tau$ function when $\phi^{(n)}$ is
as in Theorem 6.2. For $S,T\subseteq \{ 1, \dots ,N\}$, let
$K_{S,K}$ be the submatrix of $K_n$ that is indexed by $(j,k)\in S\times
T$, and let $\sharp S$ be the number of elements of $S$. 
\vskip.05in
\noindent {\bf Proposition 6.4.} {\sl  (i) Suppose that
$\phi^{(n)}:(0, \infty )\rightarrow {\bf C}$ is as in (6.20). Then}
$$\tau_{K_n} (x)=\sum_{\ell =0}^N (-1)^\ell \sum_{T, S: \sharp S=\sharp
T=\ell} \prod_{j\in S}\xi_je^{-2\lambda_jx}\prod_{k\in T}\bar\xi_k
e^{-2\bar\lambda_kx} \det\Bigl[
{{1}\over{\lambda_j+\bar\lambda_k}}\Bigr]^2_{j\in S, k\in
T}.\eqno(6.28)$$ 
\indent {\sl (ii) Suppose that $\phi^{(n)} :(0, \infty )\rightarrow E$ where
$E$ has orthonormal basis $(e_r)_{r=1}^\nu$ and let $\xi^{(r)}_j=\langle
\xi_j, e_r\rangle $. Then
$$\tau_{K_n} (x)=\sum_{S,T:\sharp S=\sharp T}(-1)^{\sharp S}\det\Bigl[
{{\xi_j^{(r)}e^{-2\lambda_jx}}\over{\lambda_j+\bar\lambda_k}}
\Bigr]_{j\in S; (k,r)\in T} \det\Bigl[
{{\bar\xi_k^{(r)}e^{-2\bar\lambda_kx}}\over{\lambda_m+\bar\lambda_k}}
\Bigr]_{m\in S; (k,r)\in T}\eqno(6.29)$$
\noindent and the sum is over all pairs of subsets 
$S\subseteq \{ 1, \dots , N\}$ and
$T\subseteq \{1, \dots , N\}\times \{1,\dots ,\nu\}$ that have
equal cardinality.}\par

\vskip.05in
\noindent {\bf Proof.} (i) By the Lemma we have $\tau_{K_n} (x)=
\det (I-K_n)$, and
by expansion of the determinant we have
$$\det (I -K_n)=\sum_{S:S\subseteq \{1, \dots , N\}} (-1)^{\sharp S}\det
K_{S,S}\eqno(6.30)$$
\noindent where $\det K_{\emptyset ,\emptyset }=1$ and otherwise
$$\det K_{S,S}=\det \Bigl[ \sum_{k=1}^N {{\xi_j\bar\xi_k
e^{-2(\lambda_j+\bar\lambda_k)x}}\over{(\lambda_j+\bar\lambda_k)(\bar
\lambda_k+\lambda_m)}}\Bigr]_{j,m\in S}\eqno(6.31)$$
\noindent which reduces by the Cauchy--Binet formula to
$$\eqalignno{\sum_{T:\sharp T=\sharp S}&\det \Bigl[
{{\xi_je^{-2\lambda_jx}}\over{\lambda_j+\bar\lambda_k}}\Bigr]_{j\in S,
k\in T}\det\Bigl[
{{\bar\xi_ke^{-2\bar\lambda_kx}}\over{\bar\lambda_k+\lambda_m}}\Bigr]_{k\in
T, m\in S}&(6.32)\cr
&=\sum_{T:\sharp T=\sharp S}\Bigl(\prod_{j\in S}\xi_je^{-2\lambda_jx}\, 
\prod_{k\in T}\bar\xi_ke^{-2\bar\lambda_k x}\Bigr) \det\Bigl[
{{1}\over{\lambda_j+\bar\lambda_k}}\Bigr]_{j\in S, k\in T}\det
\Bigl[ {{1}\over{\lambda_m+\bar\lambda_k}}\Bigr]_{m\in S, k\in
T}.\cr}$$ 
\noindent By taking the sums over both $S$ and $T$, we obtain the
stated formula.\par 
\indent (ii) To prove (ii) one follows a similar route until line
(6.32), except that we have $\langle \xi_j, \xi_k\rangle
=\sum_{r=1}^\nu\xi_j^{(r)}\bar\xi_k^{(r)}$, so the indices in the Cauchy--Binet formula are over
the product set $T\subseteq \{1, \dots , N\}\times\{ 1, \dots ,\nu\}$.
\par
\rightline{$\square$}\par

\vskip.05in



\par
\vskip.05in

\noindent {\bf 7. The $\tau$ function for the hard spectral edge}\par
\vskip.05in
\noindent Our first application of section 6 is to the hard edge
ensemble. The Jacobi polynomials arise when one applies the
Gram--Schmidt process to $(x^k)_{k=0}^\infty$ with respect to
the weight $(1-x)^\alpha (1+x)^\beta$ on $[-1,1]$ for $\alpha
, \beta >-1.$ The zeros of the polynomials of high degree tend
to accumulate at the so-called hard edges $1-$ and $(-1)+$. 
According to [28], the kernel that describes the limiting behaviour
of the joint distribution of the scaled zeros near to the hard edges 
is given by
$${{J_\nu
(2\sqrt{x})\sqrt{y}J_\nu'(2\sqrt{y})-\sqrt{x}J'_\nu(2\sqrt{x})J_\nu( 
2\sqrt{y})}\over{x-y}}=\int_0^1 J_\nu (2\sqrt{tx})J_\nu (2\sqrt{ty})\,
dt\eqno(7.1)$$
\noindent on $L^2((0,1);dt)$; here $J_\nu$ is Bessel's function of the
first kind of order $\nu$. Hence we change variables and introduce
the Hankel operators on $L^2((0,\infty ); dt)$.\par 
\vskip.05in
\noindent {\bf Proposition 7.1.} {\sl For $\nu >-1$, let 
$\phi (x)=e^{-x/2}J_\nu (2e^{-x/2})$ and let $\Gamma_\phi$ be the
Hankel integral operator on $L^2(0, \infty )$ with symbol $\phi $. 
Then Theorem 6.2 applies to
$\Gamma_{\phi}$.}\par
\vskip.05in
\noindent {\bf Proof.} From the power
series for $J_\nu$, we obtain a rapidly convergent series
$$\phi (x)=\sum_{n=0}^\infty {{(-1)^ne^{-(2n+\nu+1)x/2}}\over {n!\Gamma
(\nu +n+1)}}\qquad (x>0)\eqno(7.2)$$
\noindent giving a meromorphic transfer function
$$\hat \phi (s) =\sum_{n=0}^\infty {{(-1)^n}\over {n!\Gamma (\nu
+n+1)(s+n +(\nu +1)/2)}},\eqno(7.3)$$
\noindent for which the poles form an arithmetic progression along the
negative real axis. One can alternatively express $\hat \phi$ in
terms of Lommel's functions.\par
\indent  We choose $\lambda_n= (2n+\nu +1)/2$, so $(\lambda_n)$
gives an arithmetic progression along the positive real axis,
starting at $(\nu +1)/2>0$, and
$\sum_{n=0}^\infty\lambda_n^{-2}<\infty$. The operator $\Theta :\ell^2\rightarrow
L^2(0, \infty )$ is bounded by duality since
$$\eqalignno{\int_0^\infty \bigl\vert\sum_{n=0}^\infty
a_ne^{-\lambda_nx}\bigr\vert^2\, dx&=\sum_{n,m=0}^\infty {{a_n\bar
a_m}\over{\lambda_n+\lambda_m}}\cr
&\leq C\sum_{n=0}^\infty\vert a_n\vert^2& (7.4)\cr}$$
\noindent by Hilbert's inequality. Hence $\Gamma_\phi$ is a self-adjoint trace class
operator, and Theorem 6.2 applies.\par 
\vskip.05in
\indent We can now compute some of the finite determinants that appear in the
expansion of $\det (I-\Gamma_{\phi_{(x)}}^2)$ from Proposition 6.4.\par
\vskip.05in

\noindent {\bf Definition} {\sl (Partition).} By a partition $\lambda$ we mean a list $n_1\geq n_2\geq
\dots \geq n_\ell$ of positive integers, so that the sum $\vert
\lambda\vert =\sum_{j=1}^\ell n_j$, is split into $\ell =\ell
(\lambda )$ parts. For each $\lambda$, the symmetric
group on $\vert \lambda\vert$ letters has an irreducible unitary
representation on a complex inner product space $S_\lambda$, known as
the Specht module. For notational convenience, we introduce a null
partition $\empty$ with $\ell (\emptyset )=0$ and write
${\hbox{dim}}(S_\emptyset )=1$. \par
\vskip.05in
\noindent {\bf Proposition 7.2.} {\sl Suppose that $\nu =0$.
 Let $K=\Gamma_\phi^2$ and} $\tau
(x)={\hbox{trace}}(I-KP_{[x, \infty )})$. {\sl Then $K$ is a
trace class operator on $L^2(0, \infty )$ such that $0\leq K\leq I$ and}
$$\tau (x)=\sum_{\lambda} (-1)^{\ell (\lambda )}{{{\hbox{dim}} (S_\lambda
)^2}\over {(\vert\lambda\vert !)^2}}e^{-2\vert\lambda\vert
x}\eqno(7.5)$$
\noindent {\sl where the sum is over all partitions.}\par
\vskip.05in

\noindent {\bf Proof.} Let $E_n={\hbox{span}}\{ e^{-(2j+\nu +1)x}: j=0, \dots ,n\}$
and let $Q_n:L^2(0, \infty )\rightarrow E_n$ be the orthogonal
projection; likewise we introduce the closure $E_\infty$ of the
subspace $\cup_{n=1}^\infty E_n$ and the corresponding orthogonal
projection $Q_\infty :L^2(0, \infty  )\rightarrow E_\infty$. Observe
that $Q_n\rightarrow Q_\infty$ in the strong operator
topology as $n\rightarrow\infty$ and that 
$\Gamma_{\phi_{(x)}} Q_\infty =Q_\infty \Gamma_{\phi_{(x)}}$; hence 
$\det (I-\Gamma_{\phi_{(x)}}^2)=\lim_{n\rightarrow\infty} \det
(I-Q_n\Gamma_{\phi_{(x)}}^2Q_n)$.\par 
\indent The matrix of $Q_n\Gamma_{\phi_{(x)}}^2Q_n$ with respect to
$(e^{-(2j+\nu +1)s})_{j=0}^n$ satisfies
$$Q_n\Gamma_{\phi_{(x)}}^2Q_n\leftrightarrow\Bigl[
{{(-1)^{j+m} e^{-2x(j+m+\nu +1)}}\over {j! m!\Gamma
(\nu +j+1)\Gamma (m+\nu +1)}}\sum_{k=0}^\infty {{1}\over{(j+k+\nu
+1)(m+k+\nu +1)}}\Bigr]_{j,m =0}^{n}.\eqno(7.6)$$
\noindent We observe that the corresponding infinite matrix for $Q_\infty
\Gamma^2_{\phi_{(x)}}$ has entries that summable with respect to $j$ and
$m$ over $j,m=0, 1, \dots ;$ thus $\det (I-\Gamma^2_{\phi_{(x)}})$
is a determinant of Hill's type.\par 
\indent We consider the determinant in (6.28). We change notation so as
to allow the running indices in sums to be $j,k =0, 1, \dots, $ and we
let $S$ and $T$ be subsets of $\{ 0, 1, 2, \dots \}$ that are finite
and of equal cardinality. Suppose that the elements of
$S$ are $m_1>m_2>\dots >m_\ell $, while the elements
of $T$ are $k_1>k_2>\dots >k_\ell$; next let
$N=\ell+\sum_{i=1}^\ell (m_i+k_i)$. Then in Frobenius's
coordinates [8, 21], there is a partition $\lambda\leftrightarrow (m_1, \dots ,
m_\ell ; k_1, \dots , k_\ell )$ with $\vert\lambda\vert$ 
with a corresponding
Specht module $S_\lambda$ such that 
$$\det\Bigl[ {{1}\over{m!\Gamma (m +1)(m+k
+1)}}\Bigr]_{m\in S, k\in T}{{\prod_{k\in T}k!}\over{\prod_{m
\in S}m!}} {{{\hbox{dim}}(S_\lambda  
)}\over{(\vert\lambda\vert )!}}\eqno(7.7)$$
\noindent as in the hook length formula of representation theory; see
in [21]. Hence
the pair of sets $S$ and $T$, each with $\ell (\lambda )$ elements give
rise to the product of determinants 
$$\det\Bigl[ {{1}\over{j!\Gamma (j+1)(j+k
+1)}}\Bigr]_{j\in S, k\in T}\det\Bigl[ {{1}\over{m!\Gamma (m +1)(m+k
+1)}}\Bigr]_{m\in S, k\in T}={{{\hbox{dim}}(S_\lambda)^2}\over {(\vert
\lambda\vert !)^2}}\eqno(7.8)$$
\noindent and the exponential 
$$e^{-\sum_{j\in S}(2j+1)x-\sum_{k\in T}(2k+1)x}=
e^{-2\vert \lambda\vert x}.\eqno(7.9)$$
\indent Conversely, each partition $\lambda$ of some positive
integer gives a Ferrers diagram
and we can introduce subsets $S, T\subset \{ 0, 1, \dots \}$ that are
finite and of equal cardinality which gives a contribution to the sum
(6.28) from the prescription of (7.9) and (7.10). By summing over all partitions, or equivalently all pairs of
sets $S$ and $T$, we obtain the series (7.6).\par
\rightline{$\square$}\par
\vskip.05in 
\noindent {\bf Remark.} Borodin, Okounkov and Olshanski [8] have computed a Fredholm
determinant for the discrete Bessel kernel, and derived a result
vaguely similar to (7.6). The determinant $\det (I-KP_{(0,s)})$
was computed by Forrester, and Forrester and Witte have considered
various circular ensembles [11]. Basor and Ehrhardt have considered
asymptotics of Bessel operators [3].\par

\par
\vskip.1in

\noindent {\bf 8. A $\tau$ function related to Lam\'e's equation}\par
\vskip.05in

\indent To conclude this paper, we consider Hankel operators related to
Lam\'e's equation. First we review some ideas that originate with
Hochstadt and are developed by McKean and van Moerbecke in [23].\par 
\indent Let ${\cal E}$ be a compact Riemann surface of genus $g$, and
${\bf J}$ the Jacobi variety of ${\cal E}$, which we identify with
${\bf C}^{g}/{\bf L}$ for some lattice ${\bf L}$ in ${\bf C}^g$. An abelian function is a locally rational function on
${\bf J}$, or equivalently a periodic meromorphic function on ${\bf
C}^g$ with $2g$ complex periods. A theta function (or elliptic function
of the second kind)
$\theta :{\bf C}^g\rightarrow {\bf P}^1$ with respect to ${\bf L}$
is a meromorphic function, not identically zero, such that
there exists a linear map $x\mapsto  L(x,u)$ for $x\in {\bf C}^g$ and
$u\in {\bf L}$ and a function $\eta :{\bf L}\rightarrow {\bf C}$ such
that $\theta (x+u)=\theta (x)e^{2\pi i(L(x,u)+\eta (u))}$ for all $x\in
{\bf C}^g$ and $u\in {\bf L}$. The pair $(L, \eta )$ is called the
type of $\theta$, as in [20].\par  
\indent Suppose that $q:{\bf R}\rightarrow {\bf R}$ is infinitely
differentiable and periodic with period one. Let $U_\lambda $ be the fundamental
solution matrix for Hill's equation
$$-{{d^2}\over{dt^2}}f+q(t)f(t)=\lambda f(t)\eqno(8.1)$$
\noindent so that $U_\lambda (0)=I$, and let $\Delta (\lambda
)={\hbox{trace}}\,U_\lambda (1)$ be the discriminant. Suppose in particular that $\lambda$ lies inside the
Bloch spectrum of
$-{{d^2}\over{dt^2}}+q(t)$, but that $4-\Delta (\lambda )^2\neq 0$.
Then any nontrivial solution of (8.1) is bounded but not periodic.\par
\par
\indent We suppose that
$4-\Delta (\lambda )^2$ has only finitely many simple zeros
$0<\lambda^{(1)}_0<\lambda_1^{(1)}<\dots <\lambda_{2g}^{(1)}$, and let 
$\lambda_k^{(2)}$ be
double zeros for $k=1, 2,\dots ;$ then
$$4-\Delta (\lambda )^2=c_1\prod_{j=0}^{2g} \Bigl(
1-{{\lambda}\over{\lambda^{(1)}_j}}\Bigr)\prod_{k=1}^\infty \Bigl(
1-{{\lambda}\over{\lambda_k^{(2)}}}\Bigr)^2.\eqno(8.2)$$
\vskip.05in
\noindent {\bf Proposition 8.1.} {\sl Suppose that $q$ is a finite
gap potential so that the discriminant has this
form. Then Hill's equation gives a
Tracy--Widom system on a hyperelliptic curve of genus $g$.}\par 
\vskip.05in

\noindent {\bf Proof.} The equation (8.1) has nontrivial bounded
solutions if and only if $\vert\Delta (\lambda )\vert <2$, so that
$\lambda $ lies in an interval of stability. Hence the spectrum of $-{{d^2}\over{dt^2}}+q$ in
$L^2({\bf R})$ has
the form 
$$[\lambda^{(1)}_0, \lambda^{(1)}_1]\cup [\lambda^{(1)}_2, 
\lambda^{(1)}_3]\cup\dots
\cup [\lambda^{(1)}_{2g}, \infty ).\eqno(8.3)$$
\noindent The zeros of $\Delta'(\lambda )$ consist of all the
$\lambda_k^{(2)}$ together with zeros
$\lambda'_j$ that interlace the simple zeros of $4-\Delta (\lambda
)^2$, so $\lambda^{(1)}_{2j-1}<\lambda_j'<\lambda_{2j}^{(1)}$ for
$j=1, \dots ,g$; hence
$${{\Delta'(\lambda )}\over{\sqrt{4-\Delta (\lambda
)^2}}}={{\prod_{j=1}^g\Bigl(
1-{{\lambda}\over{\lambda'_j}}\Bigr)}
\over{\sqrt{ \prod_{j=0}^{2g}
\Bigl( 1-{{\lambda}\over{\lambda^{(1)}_j}}\Bigr)}}}.\eqno(8.4)$$
\noindent We introduce the hyperelliptic curve 
$${\cal E}: \qquad  Z^2=\prod_{j=0}^{2g}\Bigl(
1-{{X}\over{\lambda^{(1)}_j}}\Bigr),\eqno(8.5)$$
\noindent which has genus $g$. 
We introduce a new variable by the integral 
$$t=-\int {{\Delta'(X)dX}\over{\sqrt{4-\Delta (X
)^2}}}\eqno(8.6)$$
\noindent so that $2\cos t=\Delta (X)$, then we invert this relation by
introducing a hyperelliptic function $Q(t)$ with local inverse $R$ so 
that $R(Q(t))=t$ and $2\cos t=\Delta
(Q(t))$. After a little reduction, Hill's equation becomes
$$-\Bigl( {{Z}\over{\prod_{j=1}^g\bigl( 1-X/\lambda_j')
}}{{d}\over{dX}}\Bigr)^2 f+q(R(X))f=\lambda f.\eqno(8.7)$$
\noindent Now by [23, p. 260], $q(R(X))$ 
is an abelian function on ${\cal E}$ and may be viewed as a locally
rational function on the Jacobian variety ${\bf J}$ over ${\cal E}$; 
hence we can express (8.7) as a matrix differential
equation with coefficients in the field of locally rational functions on
${\bf J}$.\par
\rightline{$\square$}\par
\vskip.05in
\indent Suppose in particular that $q$ is elliptic with periods $2K$
and $2K'i$ where $K, K'>0$. Gesztesy and Weikard [12] have shown that the
spectrum has only finitely many gaps if and only if $z\mapsto
U_\lambda (z) $ is meromorphic (and possibly multivalued) for all $\lambda\in {\bf C}$. By a
classical result of Picard, there exists a nonsingular matrix
$A_\lambda$ such that $U_\lambda (z+2K)=U_\lambda (z)A_\lambda$. If
$A_\lambda$ has distinct eigenvalues, then there exists a solution $f$ to
(8.1) that is a theta function with respect to the lattice ${\bf L}=\{
2Km+2K'in:m,n\in {\bf Z}\}$. \par
\indent Next we describe in more detail the case of genus one. We recall Jacobi's sinus amplitudinus of modulus $k$ is 
${\hbox{sn}}(x\mid k)=\sin \psi$ where
$$x=\int_0^\psi {{d\theta}\over{\sqrt{1-k^2\sin^2\theta}}}.\eqno(8.8)$$
\indent For $0<k<1$, let $K (k)$ be the complete elliptic
integral
$$K(k)=\int_0^{\pi/2}{{dt}\over{\sqrt{1-k^2\sin^2t}}};\eqno(8.9)$$
\noindent next let $K'(k)=K(\sqrt{1-k^2})$; then ${\hbox{sn}}(z\mid k)^2$
has real period $K$ and complex period $2iK'$. We introduce
$$\bigl( e_1, e_2, e_3)=\Bigl({{2-k^2}\over{3}}, {{2k^2-1}\over{3}},
-{{k^2+1}\over{3}}\Bigr),\eqno(8.10)$$
\noindent and 
$$g_2={{4(k^4-k^2+1)}\over{3}}, \quad
g_3={{4(k^2-2)(2k^2-1)(k^2+1)}\over{27}};\eqno(8.11)$$
\noindent then let Weierstrass's function be 
$${\cal P}(z)=e_3+(e_1-e_2)\bigl({\hbox{sn}}(z\mid
k)\bigr)^{-2}.\eqno(8.12)$$
Likewise, ${\cal P}(z)$ has periods $2K$ and $2iK'$, and
${\cal P}(x+iK')$ is bounded, real and $2K$-periodic. In terms
of the new variable $x=z+iK'$ and the constant $B=-\lambda (e_1-e_2)-\ell
(\ell +1)e_3$, Lam\'e's differential equation (1.16) transforms to
$$\Bigl( -{{d^2}\over{dx^2}}+\ell (\ell +1) {\cal
P}(x)\Bigr)\Phi (x)+B\Phi (x)=0.\eqno(8.13)$$
\indent Writing $X={\cal P}(x)$, $Y={\cal P}(y)$ and $Z={\cal P}'(x),$ the point
$(X,Z)$ lies on the elliptic curve
$${\cal E}:\quad Z^2=4(X-e_1)(X-e_2)(X-e_3)\eqno(8.14)$$
\noindent and the elliptic function field ${\bf K}$ consists of the
field of rational functions of $X$ with $Z$ adjoined and we think of
 $B$ as a point on ${\cal E}$.
For $(x_0, z_0)$ on ${\cal E}$, we introduce the function
$$\Phi (X,Z; x_0,z_0)=\exp\Bigl( {{1}\over{2}}\int_{\gamma}
{{z-z_0}\over{x-x_0}}{{dx}\over{z}}\Bigr),\eqno(8.15)$$
\noindent which takes multiple values depending upon the path
$\gamma$ from 
$(x_0, z_0)$ to $(X,Z)$. Then for integers $\ell \geq 1$, and
typical values of $B$, there exist $\kappa\in
{\bf C}$ and polynomials $A_0(X)$ and $A_1(X)$ such that
$$\Psi (X)=\Bigl( A_0(X)+A_1(X)\Bigl( {{Z+z_0}\over{X-x_0}}\Bigr)\Bigr)\Phi
(X,Z;x_0,z_0)\exp\Bigl( \kappa \int_\gamma {{dx}\over{z}}\Bigr)\eqno(8.16)$$
\noindent gives a solution of 
$$-\Bigl( Z{{d}\over{dX}}\Bigr)^2\Psi (X)+\ell (\ell +1)X\Psi
(X)+B\Psi (X)=0,\eqno(8.17)$$
\noindent known as a Hermite--Halphen solution. Maier [22, Theorem 4.1] has
shown how to compute $(x_0,y_0)$ and the spectral curve in terms of
$\kappa$ and $B$, thus making (8.17) convenient for computation. As
$Z$ is rational on the elliptic curve, Lam\'e's equation gives 
rise to a Tracy--Widom system (1.1) that closely resembles the Laguerre
system of orthogonal polynomials with parameter one, 
as considered in [5, 29].\par
\indent Suppose henceforth that $\ell =1$. For 
$\lambda\in [k^2,1]\cup [k^2+1, \infty )$,
all solutions to (1.16) are bounded; however, except for the countable
subset of values of $\lambda$ that gives the periodic spectrum, these
solutions are not $K$ or $2K$ periodic; see [22]. Write 
$B={\cal P}(\alpha )$ where $\alpha$
is the spectral parameter. Weierstrass introduced the functions
$$\sigma (z)=z\prod_{\omega\in {\bf L}^*}\Bigl(
1-{{z}\over{\omega}}\Bigr)\exp\Bigl(
{{z}\over{\omega}}+{{1}\over{2}}\Bigl(
{{z}\over{\omega}}\Bigr)^2\Bigr)\eqno(8.18)$$
 \noindent where ${\bf L}^*={\bf L}\setminus \{ (0,0)\}$, and $\zeta (z)=\sigma'(z)/\sigma (z)$ so that ${\cal
P}=-\zeta'$. Then by [19, (13)] the equation (8.13) has a nontrivial solution 
$$\Psi (x;\alpha )=-{{\sigma (x-\alpha )}\over{\sigma (\alpha )\sigma (x)}}e^{\zeta
(\alpha )x}\eqno(8.19)$$ 
\noindent such that $\Psi (x;\alpha )\Psi
(-x;\alpha )={\cal P}(\alpha )-{\cal P}(x)$ and $\alpha \mapsto \Psi
(x;\alpha )$ is doubly periodic.\par 
\indent  The solutions
give rise to a natural kernel, for after we make the local change
of independent variable $x\mapsto X$ and write $f(X)=\Psi (x;\alpha )$ and 
$g(X)=\Psi'(x;\alpha )$, we have by [19, (18)]
$${{f(X)g(Y)-g(X)f(Y)}\over {X-Y}}=\Psi (x+y; \alpha ).\eqno(8.20)$$
\noindent The right-hand side has the shape of the kernel of Hankel
integral operator. In the remainder of this section we introduce this
operator, and compute the corresponding Fredholm determinant.
\vskip.05in
\noindent {\bf Lemma 8.2.} {\sl  Let $\beta =-2K\zeta
(\alpha )+\alpha\zeta (\alpha +2K) -\alpha \zeta (\alpha )$, suppose that $\Re \beta >0$ and 
let $t\in {\bf C}$ such that $\Psi
(x+2t;\alpha )$ is analytic for $x\in [0,2K]$. Let $\phi_{(t)} (x)
=\Psi (x+2t;\alpha )$
and $h(s)=\int_0^{2K}e^{-su}\phi_{(t)} (u)\, du$.
Then $\phi_{(t)}$ is a theta function and has an exponential expansion}
$$\phi_{(t)} (x)=\sum_{m=-\infty}^\infty {{1}\over{2K}}h\Bigl( {{2\pi
im-\beta}\over {2K}}\Bigr) e^{x(2\pi im -\beta )/(2K)}\qquad
(x>0)\eqno(8.21)$$
\noindent {\sl and $\hat \phi_{(t)}$ is a meromorphic function with poles in
an arithmetic progression.}\par
\vskip.05in
\noindent {\bf Proof.}  We introduce $\eta =\zeta (\alpha +2K)-\zeta
(\alpha )$ and
$\eta' =\zeta (\alpha +2iK')-\zeta (\alpha )$. Then $\sigma$ is a theta function and
satisfies a simple functional equation given in [20, p.109]; from this
we deduce that $\Psi$ is also a theta
function  and satisfies the functional
equations
$$\Psi (x+2K; \alpha )=\Psi (x;\alpha )e^{2K\zeta (\alpha )-\alpha\eta }, 
\quad \Psi
(x+2iK';\alpha )=\Psi (x;\alpha )e^{2iK'\zeta (\alpha )-\alpha \eta'}.
\eqno(8.22)$$
\noindent Hence $x\mapsto \Psi
(x+2t;\alpha )$ is of exponential decay as $x\rightarrow\infty$ through
real values.\par
\indent Due to (8.22), the transfer function of $\phi_{(t)} (x)$ is
$$\eqalignno{\hat \phi_{(t)} (s)&=\sum_{k=0}^\infty
\int_{2Kk}^{2K(k+1)}e^{-su}\Psi (u+2t;\alpha )\, du\cr
&=(1-e^{-2Ks+2K\zeta (\alpha )-\alpha \eta })^{-1}\int_0^{2K}\Psi
(u+2t; \alpha ) e^{-su}\, du&(8.23)}$$
\noindent which is meromorphic with possible poles at
the points $s=(2K)^{-1}(2K\zeta (\alpha )-\alpha \eta +2\pi m i)$ for $m\in {\bf
Z}$ which form a vertical arithmetic progression in the left half
plane. The position of the poles is determined by the type of the
theta function.\par
\indent We can deduce the exponential expansion by inverting the
Laplace transform. Let $T=(2m+1)\pi /(2K)$ let $x>0$ and consider the contour
$[-iT,iT]\oplus S_T$, where $S_T$ is the semicircular arc in the left
half plane with centre $0$ that goes from $-iT$ to $iT$; then by
Cauchy's Residue Theorem we have
$$\int_{S_T}e^{sx}\hat \phi_{(t)} (s)\, ds+\int_{[-iT,iT]}e^{sx}\hat \phi_{(t)}
(s)\, ds={{\pi i}\over{K}}\sum_{n=-m}^m h\Bigl({{2\pi ni-\beta
}\over{2K}}\Bigr) e^{x(2\pi ni-\beta )/(2K)}.\eqno(8.24)$$
\noindent We integrate $\int_0^{2K}\Psi (u+2t;\alpha )e^{-su}\,du$ by parts and write 
$$e^{sx}\hat \phi_{(t)} (s)={{e^{sx}}\over{s(1-e^{-2Ks-\beta})}}\Bigl(
-e^{-2Ks}\phi_{(t)} (2K)+\phi_{(t)} (0)+\int_0^{2K}e^{-su}\phi_{(t)}'(u)\,
du\Bigr)\eqno(8.25)$$
\noindent and then use Jordan's Lemma to show that
$\int_{S_T}e^{sx}\hat\phi_{(t)} (s)\, ds\rightarrow 0$ as
$T\rightarrow\infty$. Hence
$$\phi_{(t)} (x)={{1}\over{2\pi i}}\int_{-i\infty}^{i\infty }e^{sx}\hat\phi_{(t)}
(s)\, ds=\sum_{n=-\infty}^\infty {{1}\over{2K}}h\Bigl({{2\pi ni-\beta
}\over{2K}}\Bigr) e^{x(2\pi in-\beta )/(2K)}.\eqno(8.26)$$   
\rightline{$\square$}\par
\vskip.05in    
\noindent {\bf Theorem 8.3.} {\sl Let $\phi_{(t)}(x)= \Psi
(x+2t;\alpha )$ and let $\Gamma_{\phi_{(t)}}$ be the
Hankel integral operator on $L^2(0, \infty )$ with symbol
$\phi_{(t)}$. Then the conclusions of Theorem 6.1 hold for 
$\Gamma_{\phi_{(t)}}$.}\par

\vskip.05in
\noindent {\bf Proof.} Let $\lambda_n = (2\pi in+\beta
)/(2K)$ where $\Re \beta >0$. Then by a standard argument
from the calculus of residues, we have
$$\sum_{k=-\infty}^\infty {{1}\over{\vert
\lambda_j+\lambda_k\vert^{2}}}={{K^2\Re \coth \beta }\over{\Re \beta }} \qquad
(j\in {\bf Z}).\eqno(8.27)$$ 
\noindent \par 
\indent The operator $\Theta :L^2(0, \infty )\rightarrow
\ell^2$ given by
$$f\mapsto \Bigl(\int_0^\infty e^{-\bar\lambda_j s}f(s)\,
ds\Bigr)_{j=-\infty}^\infty\eqno(8.28)$$
is bounded. Indeed, we observe that the sequence 
$(e^{-\lambda_nx})_{n=-\infty}^\infty$ forms a
Riesz basic sequence in $L^2(0, \infty )$, in the sense that there
exists a constant $C>0$ such that
$$C^{-1}\sum_{n=-\infty }^\infty \vert a_n\vert^2 
\leq \int_0^\infty \Bigl\vert\sum_{n=-\infty }^\infty a_n e^{-\lambda_n
x}\Bigr\vert^2\, dx\leq C\sum_{n=-\infty }^\infty \vert a_n\vert^2
\eqno(8.29)$$
\noindent for all $(a_n)\in \ell^2$. To prove this, one uses
a simple scaling argument and orthogonality of the sequence $(e^{2\pi
inx})_{n=-\infty}^\infty $ in $L^2[0,1].$ In particular, this shows
that $\Theta^\dagger :\ell^2\rightarrow L^2(0, \infty )$ is bounded,
so $\Theta$ is bounded.\par
\indent We can now use the general Theorem 6.1. Given this rapid decay and the fact that $\Psi
(x+y+2t;a)$ is analytic, one can easily check that
$\Gamma_{\phi_{(t)}}$ is
trace class.\par
\rightline{$\square$}\par
\noindent Our final result gives the order of growth of the determinant
$$D_N=\det \Bigl[
{{1}\over{\lambda_j+\bar\lambda_k}}\Bigr]_{j,k=1}^N.\eqno(8.30)$$

\vskip.05in 
\noindent {\bf Proposition 8.4.} {\sl Suppose that $\lambda_j=(2\pi
ij+\beta )/(2K)$ where $\Re\beta >0$ and $K>0$. Let $\mu$ be the Haar 
probability measure on the
unitary group $U(N)$, and let $\arg e^{i\theta}=\theta$ for $0<\theta
<2\pi$.\par
\noindent (i) Then}
$$D_N=\Bigl(
{{2K}\over{1-e^{-2\Re\beta}}}\Bigr)^N\int_{U(N)}\exp \Bigl( -{{\Re\beta}\over{\pi}} 
{\hbox{trace}} \arg U\Bigr)\mu (dU).\eqno(8.31)$$
\indent {\sl (ii) There exists a constant $c>0$ such that}
$$\Bigl( {{K}\over
{\sinh\Re \beta}}\Bigr)^Ne^{-(2c)^{1/3}N^{2/3}(\Re \beta)^{2/3}}\leq 
D_N\leq 
\Bigl( {{K}\over
{\sinh\Re \beta}}\Bigr)^Ne^{(2c)^{1/3}N^{2/3}(\Re \beta)^{2/3}}.
\eqno(8.32)$$
\noindent {\sl so}
$$D_N^{1/N}\rightarrow K{\hbox{cosech}}\, \Re\beta \qquad
(N\rightarrow\infty ).\eqno(8.33)$$
\vskip.05in
\noindent {\bf Proof.} (i) Let 
$$f(u)={{2Ke^{-2\Re\beta u}}\over{1-e^{-2\Re\beta}}}\qquad
(0<u<1)\eqno(8.34)$$
\noindent and let the Fourier coefficients of $f$ be $a_k=\int_0^1
f(u)e^{-2\pi iku}du$, which we compute and find 
$${{1}\over{\lambda_j+\bar\lambda_k}}=a_{j-k}.\eqno(8.35)$$
\noindent Then we can use an identity due to Heine, and express the
Toeplitz determinant of $[a_{j-k}]$ as an integral
$$\det [a_{j-k}]_{j,k=1, \dots ,
N}={{1}\over{N!}}\int_{[0,1]^N}\prod_{1\leq j<k \leq N} \bigl\vert
e^{2\pi i\theta_j}-e^{2\pi i\theta_k}\bigr\vert^2\prod_{j=1}^N
f(\theta_j) \,d\theta_1\dots d\theta_N,\eqno(8.36)$$ 
\noindent which we regard as an integral over the maximal
torus in $U(N)$, and hence we convert the expression into an integral over 
the group $U(N)$, obtaining
$$\det\Bigl[
{{1}\over{\lambda_j+\bar\lambda_k}}\Bigr]_{j,k=1}^N=\int_{U(N)}
\exp\Bigl\{ {\hbox{trace}}\log f\bigl(\arg U/(2\pi )\bigr)\Bigr\} \mu (dU).
\eqno(8.37)$$
\indent (ii) Note that $\log f(\arg e^{i\theta}/(2\pi))=\log
(2K/(1-e^{-2\Re\beta})) -\Re\beta \theta /\pi$. Let $U\in U(N)$ have eigenvalues $e^{i\theta_1}, \dots ,
e^{i\theta_N}$ where $0\leq \theta_1\leq\dots\leq\theta_N\leq 2\pi$;
then the expression
$${\hbox{trace}}\arg U-\pi N=\theta_1+\dots +\theta_N-N\pi\eqno(8.38)$$
\noindent satisfies a central limit theorem, but we need to adjust the functions
slightly to accommodate the discontinuity of $\arg$. Let $g_1,g_2:{\bf R}\rightarrow {\bf R}$ be Lipschitz functions with
Lipschitz constant $L$, that are periodic with period $2\pi$, and satisfy $g_1(\theta
)\leq \theta\leq g_2(\theta )$ for $0\leq \theta <2\pi$, and 
$$\pi -{{1}\over{L}}\leq \int_0^{2\pi} g_1(\theta )\,d\theta \leq 
\int_0^{2\pi} g_2(\theta )\,d\theta
\leq \pi +{{1}\over{L}}.\eqno(8.39)$$
\indent By Szeg\"o's asymptotic formula [18], there exists a constant $c$ such that   
$$\eqalignno{\int_{U(N)} \exp \Bigl(-{{\Re\beta
}\over{\pi}}\sum_{j=1}^N\theta_j\Bigr) \mu (dU)
&\leq \int_{U(N)} \exp\Bigl( -{{\Re \beta}\over{\pi}} \sum_{j=1}^N g_1(\theta_j)\Bigr)\mu (dU)\cr
&\leq \exp \Bigl( -N\Re \beta\int_0^{2\pi} g_1( \theta){{d\theta}\over{\pi}}+c(\Re\beta
)^2L^2\Bigr);&(8.40)\cr}$$
\noindent hence we have an upper bound on $D_N$ of 
$$\Bigl( {{2K}\over{1-e^{-2\Re\beta}}}\Bigr)^N\int_{U(N)} 
\exp \Bigl(-{{\Re\beta }\over{\pi}}\sum_{j=1}^N\theta_j\Bigr) \mu (dU)\leq
\Bigl( {{2K}\over{e^{\Re\beta
}-e^{-\Re\beta}}}\Bigr)^N e^{\Re \beta N/L +c(\Re\beta
)^2L^2}.\eqno(8.41)$$
\noindent Using $g_2$ instead of $g_1$, one can likewise obtain a lower bound on
$D_N$. To conclude the proof, we choose $L=N^{1/3}(2c\Re \beta )^{-1/3}$.\par
\rightline{$\square$}\par
\vskip.05in
\noindent {\bf References}\par
\noindent [1] M.J. Ablowitz and A.S. Fokas, {Complex Analysis:
Introduction and Applications}, 2nd Edition, (Cambridge, 2003).\par
\noindent [2] M.J. Ablowitz and I.A. Segur, Exact linearization of a
Painlev\'e transcendent, {\sl Phys. Rev. Lett.} {\bf 38} (1977),
1103--1106.\par  
\noindent [3] E.L. Basor and T. Ehrhardt, Asymptotics of determinants of
Bessel operators, {Commun. Math. Physics} {234} (2003),
491--516.\par
\noindent [4] G. Blower, Operators associated with soft and hard
spectral edges
from unitary ensembles, {J. Math. Anal. Appl.} {337} (2008),
239--265.\par
\noindent [5] G. Blower, Integrable operators and the squares of Hankel
operators, {J. Math. Anal. Appl.} {340} (2008), 943--953.\par
\noindent [6] G. Blower, Linear systems and determinantal random point
fields, {J. Math. Anal. Appl.} {355} (2009), 311--334.\par

\noindent [7] A. Borodin and P. Deift, Fredholm determinants,
Jimbo--Miwa--Ueno $\tau$--functions and representation theory, {
Comm. Pure Appl. Math.} {55} (2002), 1160--1230.\par
\noindent [8] A. Borodin, A. Okounkov and G. Olshanski, Asymptotics of
Plancherel measures for symmetric groups, {J. Amer. Math. Soc.} {
13} (2000), 481--515.\par 
\noindent [9] P.A. Deift, A.R. Its, and X. Zhou, A Riemann--Hilbert
approach to asymptotic problems arising in the theory of random matrix
models, and  also in the theory of integrable statistical mechanics,
{Annals of Math.} {(2)} {146} (1997), 149--235.\par 
\noindent [10] A.S. Fokas, A.R. Its, A.A. Kapaev and V.Y. Novokshenov,
{\sl Painlev\'e transcendents: the Riemann--Hilbert approach},
Mathematical Surveys and Monographs 128,  American
Mathematical Society 2006.\par
\noindent [11] P.J. Forrester and N.S. Witte, Applications of the
$\tau$-function theory of Painlev\'e equations to random matrices:
$P_V$, $P_{III}$, the LUE, JUE and CUE, {Comm. Pure Appl. Math.}
{55} (2002), 679--727.\par
\noindent [12] F. Gesztesy and T. Weikard, Picard's equation and Hill's
equation on a torus, {Acta Math.} {176} (1996), 73--107.\par 
\noindent [13] D. Guzzetti, The elliptic representation of the 
general Painlev\'e VI equation, {Comm. Pure Appl. Math.} {55} (2002),
1280--1363.\par
\noindent [14] M. Jimbo, Monodromy problem and the boundary condition
for some Painlev\'e equations, {Publ. Res. Inst. Math. Sci.} 18 (1982),
1137--1161.\par
\noindent [15] M. Jimbo, T. Miwa and K. Ueno, Monodromy preserving
deformations of linear ordinary differential equations with
rational coefficients I: general theory, {Physica D} {2}
(1981), 306--352.\par
\noindent [16] M. Jimbo and T. Miwa, Monodromy preserving deformations of 
linear ordinary differential equations with rational coefficients II,
{Physica D} {2}, 406--448.\par
\noindent [17] M. Jimbo and T. Miwa, Monodromy preserving deformations
of linear differential equations with rational coefficients III, 
{Physica D} {4} (1981/2), 26--46.\par
\noindent [18] K. Johansson, On Szeg\"o's asymptotic formula and
Toeplitz determinants and generalizations, {Bull. Sci. Math.} (2)
{112} (1988), 257--304.\par
\noindent [19] I.M. Krichever, Elliptic solutions of the
Kadomcev--Petviasvili equations, and integrable systems of particles, 
{Functional Anal. Appl.} {14} (1980), 282--290.\par 
\noindent [20] S. Lang, {Introduction to Algebraic and Abelian
Functions}, Second Edition, Springer--Verlag, 1982.\par
\noindent [21] I.G. MacDonald, {Symmetric functions and Hall
polynomials}, Oxford University Press, Second Edition, Clarendon Press, 
1995.\par
\noindent [22] R.S. Maier, Lam\'e polynomials, hyperelliptic
reductions and Lam\'e band structure, {Philos. Trans. R. Soc. A Math.
Phys. Eng. Sci.}
{336} (2008), 1115--1153.\par
\noindent [23] H.P. McKean and P. van Moerbeke, {The spectrum of Hill's
equation,} {Invent. Math.} {30} (1975), 217--274.\par 
\noindent [24] K. Okamoto, On the $\tau$-functions of Painlev\'e
equations, {Physica D} {2} (1981), 525--535.\par
\noindent [25] F.W.J. Olver, {Asymptotics and Special Functions},
Academic Press, New York, 1974.\par
\noindent [26] V.V. Peller, {Hankel Operators and Their
Applications,} Springer, New York, 2003.\par
\noindent [27] T. Stoyanova, Non-integrability of 
Painlev\'e VI equations in the
Liouville sense, {Nonlinearity} {22} (2009), 2201--2230.\par 
\noindent [28] C.A. Tracy and H. Widom, Level spacing
distributions and the Bessel kernel, {Commun. Math. Phys.}
{161} (1994), 289--309.\par 

\noindent [29] C.A. Tracy and H. Widom, Fredholm determinants,
differential equations and matrix models, {Commun. Math. Phys.}
{163} (1994), 33--72.\par 

\noindent [30] C.A. Tracy and H. Widom, Fredholm determinants and the
mKdV/sinh-Gordon hierarchies, {Comm. Math. Phys.} {179} (1996),
1--9.\par 
\noindent [31] H.L. Turrittin, Reduction of ordinary differential
equations to the Birkhoff canonical form, {Trans. Amer. Math. Soc.}
{107} (1963), 485--507.\par
\noindent [32] E.T. Whittaker and G.N. Watson, {A Course of Modern
Analysis}, fourth edition, Cambridge University Press, 1965.\par
\vfill
\eject
\end